\documentclass[11pt,a4paper]{amsart} 

\usepackage[utf8]{inputenc}
\usepackage{xcolor}
\usepackage{graphicx}
\usepackage{parskip}
\usepackage{amssymb, amsthm,amsmath}
\calclayout
\usepackage[T1]{fontenc}
\usepackage[english]{babel}
\usepackage{amsmath, amssymb, mathtools, amsthm}
\usepackage{mathrsfs}
\usepackage{dsfont}
\usepackage{stmaryrd}
\usepackage{setspace} 
\usepackage{float}
\usepackage{cleveref}
\usepackage{csquotes}
\usepackage[style=alphabetic]{biblatex}
\usepackage{extarrows}
\usepackage{longtable}
\usepackage{lineno}
\usepackage{lscape}
\usepackage{romanbar}
\usepackage{multicol}
\usepackage{wasysym}
\usepackage{stmaryrd}
\usepackage{wrapfig}
\usepackage{subcaption}
\usepackage{xfrac}
\usepackage{faktor}
\usepackage[obeyFinal]{todonotes}
\usepackage{overpic}
\usepackage{tikz}
\usetikzlibrary{positioning}
\usetikzlibrary{fadings}
\usetikzlibrary{shapes.symbols}
\usetikzlibrary{shapes.geometric}
\usetikzlibrary{decorations}
\usetikzlibrary{decorations.shapes}
\usetikzlibrary{decorations.pathmorphing}
\usetikzlibrary{arrows.meta}

\definecolor{darkturquoise}{RGB}{46, 114, 132}
\definecolor{lightturquoise}{RGB}{53, 130, 151}
\definecolor{pastelblue}{RGB}{176, 213, 226}

\usepackage{algorithm}
\usepackage{algpseudocode}

\theoremstyle{definition}
\newtheorem{model}{Definition}[section]
\theoremstyle{plain}
\newtheorem{lemma}[model]{Lemma}

\newtheorem{theorem}[model]{Theorem}
\newtheorem*{theorem*}{Theorem}
\newtheorem{proposition}[model]{Proposition}
\newtheorem{cor}[model]{Corollary}

\theoremstyle{definition}
\newtheorem{defi}[model]{Definition}

\newtheorem{example}[model]{Example}

\theoremstyle{remark}
\newtheorem{remark}[model]{Remark}

\newcommand{\id}{\mathds{1}}
\newcommand{\N}{\mathbb{N}}
\newcommand{\R}{\mathbb{R}}
\newcommand{\Z}{\mathbb{Z}}

\newcommand{\labs}{\left|}
\newcommand{\rabs}{\right|}

\makeatletter
\newcommand{\oset}[3][0ex]{%
	\mathrel{\mathop{#3}\limits^{
			\vbox to#1{\kern-2\ex@
				\hbox{$\scriptstyle#2$}\vss}}}}
\makeatother

\def\card{\operatorname{card}}
\def\sgn{\operatorname{sgn}}  
\def\d{\operatorname{d}}
\newcommand{\cf}{\mathbf{c}_\mathbf{f}}
\newcommand{\wall}{{H}}
\newcommand{\setwalls}{\mathcal{H}}
\newcommand{\cone}{\mathcal{C}}
\newcommand{\graph}{\mathcal{G}}
\newcommand{\length}{\ell}
\newcommand{\Cox}{{W}}
\newcommand{\SCox}{S}
\newcommand{\sW}{W_{0}}
\newcommand{\aW}{W}
\newcommand{\tA}{\mathtt{A}}
\newcommand{\tB}{\mathtt{B}}
\newcommand{\atA}{\tilde{\mathtt{A}}}
\newcommand{\atB}{\widetilde{\mathtt{B}}}
\newcommand{\mb}{\mathbf}

\title{Folded galleries and moment graphs}
\author{Anna Reimann}
\address{Anna Reimann\\ 
Otto von Guericke University Magdeburg \\ 
Faculty of Mathematics, Universit\"atsplatz 2, 39106 Magdeburg, Germany}
\email{anna.michael@ovgu.de}
\date{\today}

\begin{document}

\begin{abstract}
 We characterize \emph{folding patterns}, the combinatorial options of folding minimal alcove-to-alcove galleries in affine Coxeter complexes positively with respect to Weyl chamber orientations of the Coxeter complex, by drawing a connection to the Bruhat moment graph of the associated spherical Coxeter group.
 We also prove how to determine the spherical direction of the end alcove of a positively folded gallery using these graphs.
\end{abstract}

\maketitle

\section{Introduction}\label{sec:introduction}

\begin{figure}
    \begin{subfigure}[c]{0.45\textwidth}
        \centering
     \begin{overpic}[scale=0.4]{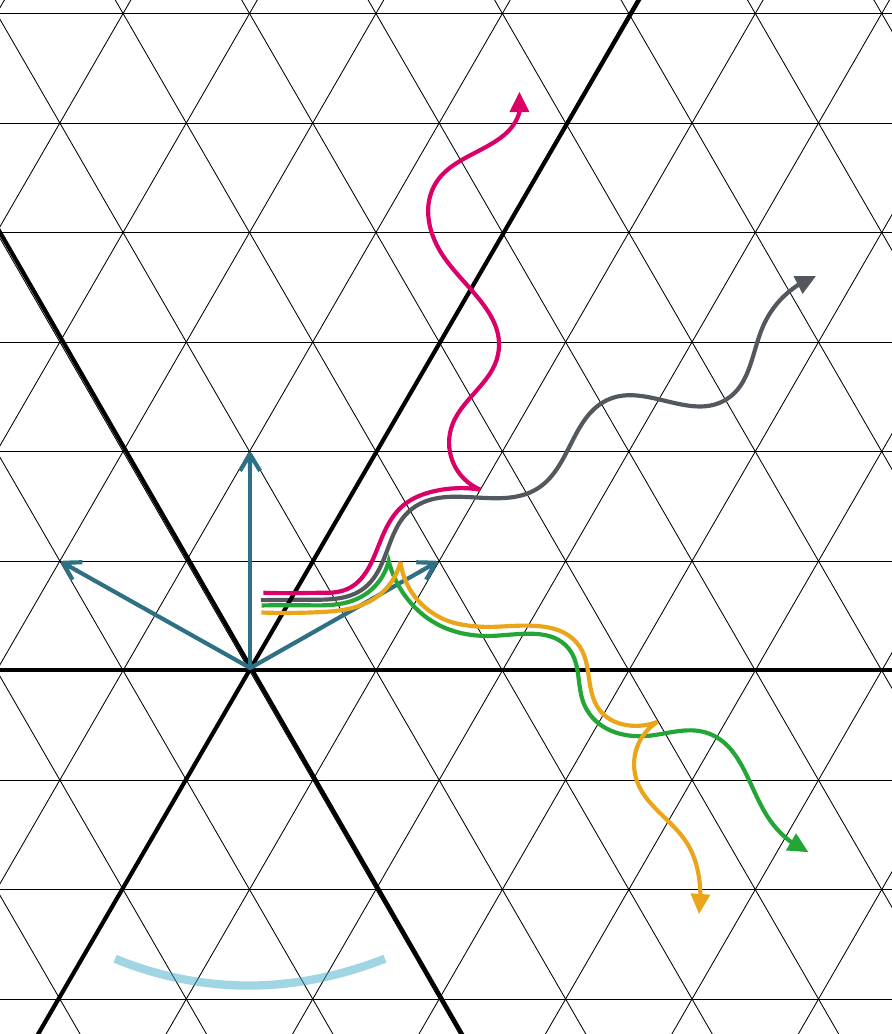}
        \put(9,39) {\tiny \textcolor{darkturquoise}{$\alpha_2$}}
        \put(34,39) {\tiny \textcolor{darkturquoise}{$\alpha_1$}}
        \put(10,49) {\tiny \textcolor{darkturquoise}{$\alpha_1+\alpha_2$}}

        \put(21,42) {\tiny $e$}
        \put(17,40) {\tiny $s_1$}
        \put(14,32) {\tiny $s_1s_2$}
        \put(19,25) {\tiny $s_1s_2s_1$}
        \put(27,32) {\tiny $s_2s_1$}
        \put(28,36) {\tiny $s_2$}

        \put(71,70) {\footnotesize $\gamma$}
        \put(74,22) {\footnotesize $\gamma_1$}
        \put(62,17) {\footnotesize $\gamma_2$}
        \put(44,89) {\footnotesize $\gamma_3$}
         %orientation
        \put(80,32) {\tiny $+$}
        \put(80,36) {\tiny $-$}
        \put(15,17) {\tiny $+$}
        \put(11,18) {\tiny $-$}
        \put(8,64) {\tiny $-$}
        \put(4,63) {\tiny $+$}
     \end{overpic}
    \end{subfigure}
    \begin{subfigure}[c]{0.45\textwidth}
        \centering
    \begin{tikzpicture}[scale=1,->, >=Stealth, shorten >=1pt, shorten <=1pt, thick]   
        %nodes    
        \node at (5, -0.7) {$e$};
        \draw[fill=black] (5,-1) circle (2pt);
        \node at (3.2, -2) {$s_1$};
        \draw[fill=black] (3.5,-2) circle (2pt);
        \node at (6.8, -2) {$s_2$};
        \draw[fill=black] (6.5,-2) circle (2pt);
        \node at (3, -4) {$s_1s_2$};
        \draw[fill=black] (3.5,-4) circle (2pt);
        \node at (7, -4) {$s_2s_1$};
        \draw[fill=black] (6.5,-4) circle (2pt);
        \node at (5, -5.3) {$s_1s_2s_1=s_2s_1s_2$};
        \draw[fill=black] (5,-5) circle (2pt);
        %edges
        \draw[->] (5,-1)--(3.5,-2);
        \draw[->] (5,-1)--(6.5,-2);
        \draw[->] (5,-1)--(5,-5);
        \draw[->] (3.5,-2)--(3.5,-4);
        \draw[->] (3.5,-2)--(6.5,-4);
        \draw[->] (3.5,-4)--(5,-5);
        \draw[->] (6.5,-2)--(6.5,-4);
        \draw[->] (6.5,-2)--(3.5,-4);
        \draw[->] (6.5,-4)--(5,-5);
        %roots
        \node[color=darkturquoise] at (4, -1.4) {\tiny$\alpha_1$};
        \node[color=darkturquoise] at (6, -1.4) {\tiny$\alpha_2$};
        \node[color=darkturquoise] at (4, -2.6) {\tiny$\alpha_2$};
        \node[color=darkturquoise] at (6, -2.6) {\tiny$\alpha_1$};
        \node[color=darkturquoise] at (2.95, -3) {\tiny$\alpha_1+\alpha_2$};
        \node[color=darkturquoise] at (7.05, -3) {\tiny$\alpha_1+ \alpha_2$};
        \node[color=darkturquoise] at (4, -4.6) {\tiny$\alpha_2$};
        \node[color=darkturquoise] at (6, -4.6) {\tiny$\alpha_1$};
        \node[color=darkturquoise] at (5.55, -4) {\tiny $\alpha_1+\alpha_2$};
    \end{tikzpicture}
    \end{subfigure}
    \caption{\emph{On the left:} Coxeter complex of $\atA_2$ with folded galleries. \emph{On the right:} Bruhat moment graph of $\tA_2$.}
    \label{fig:BruhatMomentGraphAndPosFolding}
\end{figure}

Folded galleries were introduced first by Peter Littelmann: 
He used them in his \emph{path model} to answer representation theoretic questions on symmetrizable Kac-Moody algebras \cite{Littelmann1994,Littelmann1995}.
He used the dimension of galleries to compute weight multiplicities of irreducible highest weight representations and determined the irreducible components of the tensor product of two highest weight representations of semisimple complex Lie algebras.
Together with Gaussent \cite{GaussentLittelmann_2005} they proved correspondencies between folded galleries and MV-Polytopes.
They interpreted folded galleries in the Bruhat-Tits buildings associated with a semisimple algebraic group and connected them to the structure of the corresponding affine Grassmannian.
Since then, folded galleries have shown to be a valuable combinatorial tool for studying problems of broadly diversified topics, compare \cite{Ram_2006,Kapovich2004,Schwer_2006,Ehrig2010,MilicevicSchwerThomas_2015}.\\
It then was a natural question to ask for the set of (end alcoves of) galleries that can be obtained from a given (minimal) gallery $\gamma$ by folding it.
This was approached by Graeber and Schwer in \cite{Graeber_2020}, starting to study folded galleries as combinatorial objects on their own, without pursuing a case of application, and introducing the notion of a \emph{Coxeter shadow} for this set.
So far, these Coxeter shadows can be computed using recursive algorithms, but they lack of a closed description.\\
It has been an open question to determine at which panels a minimal gallery can be folded at positively with only using the information given by the corresponding word and group element.
We answer this question for a certain type of orientation of the Coxeter complex, using the Bruhat moment graph of the associated spherical Weyl group $\sW$.
To the authors knowledge, moment graphs were first introduced to study the topology of a complex equivariantly formal variety by Goresky, Kottwitz and MacPherson \cite{GoreskyKottwitzMacPherson1997}, and later used by Peter Fiebig to interpret Lusztig's and Kazhdan-Lusztig's conjecture as problems on moment graphs \cite{Fiebig2007,Fiebig2010}, compare also the survey \cite{Fiebig2016}.\\
To obtain our main result, we use moment graphs associated with root systems of finite Coxeter groups, so called Bruhat graphs.
Notice that the Bruhat graph as used in \cite{BjornerBrenti_2006} also illustrates a poset structure on the elements of the Coxeter group, defined by the reflections of the group.
For this text, we will call a gallery $\gamma= (\mb{c}_0, p_1, \mb{c}_1, \dots, p_n, \mb{c}_n)$ with alcoves $\mb{c}_i$ and panels $p_i$ \emph{folded at $p_i$}, if $\mb{c}_{i-1}=\mb{c}_i$.
See \Cref{fig:BruhatMomentGraphAndPosFolding} for an example:
The orange gallery $\gamma_2$ is folded at its third and ninth panel.
We call the sequence of roots perpendicular to the hyperplanes containing the panels a gallery is folded in the \emph{folding pattern}.
Hence the orange gallery $\gamma_2$ has the folding pattern $(\alpha_1 + \alpha_2, \alpha_1)$.
Using this, and the correspondence of spherical and affine Coxeter groups, we will prove as \Cref{thm:positiveFoldingPatternsAndModMomentGraph} below:

\begin{theorem*}
    Fix an affine Coxeter system $(\aW, \SCox)$ with associated spherical group $\sW$.
    Let the affine Coxeter complex be equipped with a Weyl chamber orientation $\phi_{w}$, $w \in \sW$ and let $\gamma$ be a minimal gallery starting in $\cf$.
    Then the directed paths in the Bruhat moment graph of $\sW$ determine the possible $\phi_w$-positive sequences of folding hyperplanes of $\gamma$.
\end{theorem*}

An example is shown in \Cref{fig:BruhatMomentGraphAndPosFolding}:
The green and yellow galleries $\gamma_1$, $\gamma_2$ can be obtained from the unfolded minimal gallery $\gamma$ by folding positively with respect to the Weyl chamber orientation defined by $s_1s_2s_1$, marked in blue.
The galleries are folded in hyperplanes perpendicular to the roots $\alpha_1 + \alpha_2$, and $\alpha_1 + \alpha_2, \alpha_1$, respectively.
This coincides with the labels of a directed path in the Bruhat moment graph, depicted on the right in \Cref{fig:BruhatMomentGraphAndPosFolding}, starting in the node labeled $s_2$, which is the Weyl chamber $\gamma$ ends in.
As well one can see that $\gamma_3$ is not positively folded with respect to this orientation, since is is folded in a hyperplane perpendicular to the root $\alpha_2$, and there is no directed edge leaving the node $s_2$ in the graph.

Another natural question would be to ask for inclusion of Coxeter shadows: Given a group element $x$, is this element contained in the Coxeter shadow of another element $y$? That is, can a minimal gallery ending in $\bf{y}$ be positively folded with respect to a given orientation, such that the folded gallery ends in $\bf{x}$?
Regarding this question, an interesting observation on folded galleries is that folding obviously affects the spherical direction of the end alcove of the gallery.
We show a connection between the Bruhat moment graph of the associated spherical Coxeter group and the spherical direction of end alcoves of galleries in \Cref{prop:UndirectedMomentGraphAndSphericalDirection},
allowing us to answer these questions easily.

\textbf{This paper is organized as follows:} \Cref{sec:CoxSystemComplOrient} sets our notation on Coxeter systems, Coxeter complexes, and their orientation.
In \Cref{sec:PosFoldedGalleries} we state our specific definition of galleries and folds in Coxeter complexes, introduce the notion of a folding pattern, and prove technical lemmata on hyperplane crossings of these galleries needed later in the text.
Eventually, \Cref{sec:momentGraphs} points our two connections between moment graphs and folded galleries, also containing full statement and prove of the Theorem mentioned above.

\textbf{Aknowledgements:} This work was done with full financial support by the DFG, German Research Foundation, grant no. 314838170, GRK 2297 MathCoRe.
The author thanks Petra Schwer for fruitful discussions, and Jacinta Torres for helpful comments on an earlier version of this text.
%
%%%%%%%%%%%%%%%%%%%%%%%%%%%%%%%%%%%%%%%%%%%%%%%%%%%%%%%%%%%%%%
%%%%%%%%%%%% SECTION 2 %%%%%%%%%%%%%%%%%%%%%%%%%%%%%%%%%%%%%%%
%%%%%%%%%%%%%%%%%%%%%%%%%%%%%%%%%%%%%%%%%%%%%%%%%%%%%%%%%%%%%%%
%
\section{Coxeter systems and oriented Coxeter complexes}\label{sec:CoxSystemComplOrient}

In this section, we will recall basic definitions and properties, and fix notation.
For more background on Coxeter groups and their geometry, the reader may refer to one of the many good textbooks on the topic; the author especially recommends \cite{Davis_2012}, \cite{AbramenkoBrown_2010}, \cite{Humphreys_1990} and \cite{BjornerBrenti_2006}.
%
%%%%%
%
\subsection{Coxeter systems and complexes}\label{subsec:CoxSysCompl}

Let $(\Cox,\SCox)$ be a Coxeter system of rank $n+1$.
Denote by $\Sigma = \Sigma (\Cox,\SCox)$ the geometric realization of the simplicial complex associated with $(\Cox,\SCox)$.
In case $\Cox$ is affine, we can identify $\Sigma$ with a tiled $n$-dimensional Euclidean vector space the group $\Cox$ acts on.
We will call the hyperplanes $\wall$ of this tiling \emph{walls}.
$\Cox$ contains the \emph{reflections} $ r \in R:= \bigcup_{x \in \Cox} x \SCox x^{-1}$, each of them fixing a wall $\wall_r$ in $\Sigma$ which we call the \emph{reflection hyperplane of $r$}.\\
The elements of the \emph{Coxeter complex} $\Sigma(\Cox,\SCox)$ are the cosets  $x\Cox_{\SCox '}$ of parabolic subgroups of $(\Cox,\SCox)$ generated by $\SCox ' \subseteq \SCox$.
They form a poset, ordered by reverse inclusion.
The $0$-dimensional elements of $\Sigma$, the \emph{vertices}, correspond to the cosets of maximal parabolic subgroups to generating sets $\SCox ' = \SCox \setminus \lbrace s \rbrace$.
We will call the maximal dimensional simplices of $\Sigma$ \emph{alcoves} and denote them by bold lowercase letters, e.g. $\mathbf{c}$ or $\mathbf{x}$.
Observe that they are in 1-to-1 correspondence with the group elements of $\Cox$, denoted as $c$ or $x$.
Their codimension one faces will be called \emph{panels}, denoted by $p$, $q$.
Since every panel corresponds to a parabolic subgroup of the form $x\Cox_{s}$ for some $s \in \SCox$, we can assign a \emph{type} $\tau (p)=s$ to every panel.

For most of this text, we will restrict to the case of affine Coxeter systems.
If $(\aW,\SCox)$ is an affine Coxeter system, the group $\aW$ has the structure of a semi-direct product $\aW= T \rtimes \sW$, with a translation group $T \cong \Z^n$ corresponding to the coroot lattice, and an associated spherical Weyl group $\sW$.
This allows us to express a group element $x \in \aW$ as $x= t^\lambda v$ with $\lambda \in T$ and $ v \in \sW$, where we call $v$ the \emph{spherical direction} of $x$.
Let $\Sigma = \Sigma (W,S)$ be the tiled $n$-dimensional Euclidean vector space $V$ that $\aW$ acts on in this setting.
We call the vertices of $\Sigma$ whose stabilizer in $\aW$ is isomorphic to $\sW$ the \emph{special vertices}.
Denote the origin of the vector space $V$ by $v_0$ and choose the geometric realization of $\Sigma$ in a way, such that $v_0$ is a special vertex, fixed by all elements in $\SCox_0 = \lbrace s_1, \dots, s_n \rbrace$, the subset of $\SCox= \lbrace s_0, s_1, \dots, s_n \rbrace$ generating the spherical Weyl group $\sW$.
The set of walls $\setwalls_{\lambda}$ through a special vertex $\lambda$ splits $V$ into $\card (W_0)$ simplicial cones.
We call the closures of the connected components $V \setminus \setwalls_{\lambda}$ \emph{Weyl chambers} and denote the Weyl chambers containing the origin $v_0$ by $\cone_v$, indexed by the spherical group element $v$ with $t^0v \subseteq \cone_v$.\\
The set of parallel classes of rays in $\Sigma$ form a sphere at the boundary of $\Sigma$ which we call the \emph{boundary sphere} $\partial \Sigma$.
This sphere is also a Coxeter complex, associated with the spherical Weyl group $\sW$, inheriting its structure from $\Sigma$ by taking the parallel classes of rays as hyperplanes in $\partial \Sigma$.
It follows, that the alcoves in $\partial \Sigma$ are precisely the parallel classes of the Weyl chambers in $\Sigma$.
We denote the Weyl chamber with basepoint at the origin $v_0$ corresponding to the parallel class of Weyl chambers that is represented by the identity element $\id \in \sW$  by $\cone_f$, and call it the \emph{fundamental Weyl chamber}.
The unique alcove of $\Sigma$ contained in $\cone_f$ having $v_0$ as a vertex will be called \emph{fundamental alcove} and denoted by $\cf$.
Then, the $\aW$ action on $\Sigma$ identifies alcoves in $\Sigma$ with group elements $w \in \aW$ as $\mb{w}= w\cf$.
Compare \Cref{fig:BruhatMomentGraphAndPosFolding}, where the six alcoves adjacent to the origin are labeled with their corresponding group elements.\\
Besides this combinatoric view on Coxeter groups and complexes, one can also take a more algebraic approach:
Let $\Phi = \lbrace \alpha_1, \dots, \alpha_n \rbrace$ be a root system in $V$ with basis $B$, such that $(W_0, S_0)$ is the associated Weyl group.
We denote the set of positive (resp. negative) roots by $\Phi^+$ (resp. $\Phi^-$).
For this text, we want to restrict to the case of irreducible root systems.
The coroot lattice will be denoted by $R^\vee = \oplus_{i=1}^n \Z \alpha_i^\vee \subseteq \Z^n$ with coroots $\alpha^\vee = \frac{2\alpha}{\langle \alpha , \alpha \rangle}$ and basis $B^\vee$, the fundamental weights and co-weights by $\omega_\alpha$ and $\omega_\alpha^\vee$.\\
The walls $\wall_{v_0}$ of the Coxeter complex $\Sigma$ through the origin $v_0$, that is, the set of reflection hyperplanes for reflections in $\Cox$ that fix this vertex, are in fact the linear hyperplanes perpendicular to the $\alpha_i \in \Phi$.
Denote the corresponding reflections across these linear hyperplanes by $s_{\alpha_i}$.
More generally, we can consider the affine hyperplanes $\wall_{\alpha, k}:= \lbrace v \in V \vert \langle \alpha , v \rangle = k \rbrace$ for each root $\alpha \in \Phi$, and the two associated half-spaces $\wall_{\alpha, k}^+ = \lbrace x \in V \vert \  \langle x, \alpha^\vee \rangle \geq k \rbrace$ and $\wall_{\alpha, k}^- = \lbrace x \in V \vert \ \langle x, \alpha^\vee \rangle \leq k \rbrace$.
We denote by $r_{\alpha,k}$ the reflection along the hyperplane $\wall_{\alpha,k}$ and abbreviate notation for the fundamental reflections to $r_{\alpha,0}= r_{\alpha}$.
Sometimes, we also index the associated generator of $\SCox_0$ (resp. $\SCox$) by the root perpendicular to their fundamental reflection hyperplane: $\SCox_0=\lbrace s_{\alpha,0}= s_{\alpha_i} \vert \  \alpha_i \in \Phi \rbrace$, $S= \lbrace s_{\alpha,0}= s_{\alpha_i} \vert \ \alpha_i \in \Phi \rbrace \cup \lbrace s_{\tilde{\alpha},1} \rbrace$.
%%%

\subsection{Geometry of affine Coxeter complexes}\label{subsec:GeomAffCoxCompl}

Because of the semidirect product structure of an affine Coxeter Group $\aW= \sW \rtimes T$ with $\sW$ a spherical Coxeter group and $T \cong \Z^n$, we can find the tiling of $\R^n$ given by the Weyl chambers with cone point $v_0$ also at every other coroot lattice point $\mu \in R^\vee$.
This gives the Coxeter complex $\Sigma$ a structure similar to infinitely many local copies of $\Sigma(\sW, \SCox_0)$ at every $\mu \in R^\vee$.
Since we want to make use of this local-to-global-geometric structure below, we introduce the following:

\begin{defi}[Spherical direction]
    Let $(\aW,\SCox)$ be an affine Coxeter system with $\sW$ the corresponding spherical Coxeter group and let $x = t^\lambda w \in \aW$.
    Then, we denote by $\zeta : \aW \to \sW$, $x \mapsto \zeta (x) = w$ the map that records the \emph{spherical direction} of the group element $x$.
\end{defi}

Also, the following notation will help us to describe our observations in affine Coxeter complexes below more easily:

\begin{defi}[Parallelism class of a hyperplane]
    Let $\wall_{\alpha, k}$ be a hyperplane in the Coxeter complex $\Sigma$ perpendicular to a positive, not necessarily simple root $\alpha \in \Phi^+$.
    We call $\alpha$ the \emph{parallelism} or \emph{root class of $\wall_{\alpha,k}$}.
    Sometimes, we also call $\wall_{\alpha,k}$ an $\alpha$-hyperplane.
\end{defi}

\begin{defi}[$i$-strip]\label{Def:i-strip}
    Let $\Sigma = \Sigma (\aW,\SCox)$ be the affine Coxeter complex corresponding to the Coxeter system $(\aW,\SCox)$ and let $\alpha \in \Phi^+$ be a positive root of $\aW$.
    Denote by $\wall_{\alpha, 0}^{\cf}$ the half-space defined by $\wall_{\alpha,0}$ that contains the fundamental alcove $\cf$ and by $\wall_{\alpha, i}^{\cf}$ for $i \in \Z$ the half-space corresponding to $\wall_{\alpha,i}$ that satisfies $\wall_{\alpha, i}^{\cf} \subseteq \wall_{\alpha, 0}^{\cf}$ or $\wall_{\alpha, i}^{\cf} \supseteq \wall_{\alpha, 0}^{\cf}$.
    Denote the second half-space defined by $\wall_{\alpha,i}$ by $\wall_{\alpha,i}^{\neg \cf}$.
    Then, we define the \emph{$i$-strip with respect to $\alpha$} to be 
    \begin{align*}
        \mathcal{S}_{\alpha,i} := \wall_{\alpha, i}^{\cf} \cap \wall_{\alpha, i+1}^{\neg \cf} .
    \end{align*}
\end{defi}

\begin{defi}[Local Weyl chamber]
    Let $(\aW,\SCox)$ be an affine Coxeter system with $(\sW, \SCox_0)$ the corresponding spherical Coxeter system.
    Let $x = t^\mu w \in \aW$ with $\mu \in R^\vee$, $w \in \sW$ and define $h_\alpha (\mu) = \langle \mu, \alpha^\vee \rangle $ for all positive roots $\alpha \in \Phi^+$.
    Denote by $\wall_{\alpha, h_\alpha(\mu)}^{\bf{x}}$ the half-space determined by the hyperplane $\wall_{\alpha, h_\alpha(\mu)}$ that contains the alcove $\mathbf{x}$.
    Then we call the cone
    \begin{align*}
        \mathcal{C}_{\mu, w} := \bigcap_{\alpha \in \Phi^+} \wall_{\alpha, h_\alpha(\mu)}^{\bf{x}}
    \end{align*} 
    the \emph{local Weyl chamber at $\mu$ with direction $w$}.
    Notice that $\mathbf{x}$ is the alcove at the tip of the cone, and $\mu$ is the cone point.
\end{defi}

\subsection{Orientations}\label{subsec:Orientations}

In this subsection, we introduce orientations of Coxeter complexes, using the Definitions of \cite[Section 3]{Graeber_2020}, and later restrict to the class used in the text below.

\begin{defi}{(Orientation of $\Sigma$)}
	An \emph{orientation} $\phi$ of $\Sigma$ is a map that assigns to every pair of alcove $\mb{c}$ and panel $p$ contained in $\mb{c}$ either $+1$ or $-1$.\\
	We say that an alcove $\mb{c}$ is \emph{on the positive (resp. negative) side of} $p$ if $\phi (p,\mathbf{c} )= +1$ (resp. $\phi (p,\mathbf{c} ) = -1$).\\
	An orientation $\phi$ is called \emph{wall-consistent} if for any wall $\wall$ in $\Sigma$ and alcoves $\mb{c},\mb{d}$ contained in the same half-space defined by $\wall$ and having panels $p,q$ in $\wall$, one has:
	\begin{align*}
		\phi (p,\mathbf{c}) = +1 \Leftrightarrow \phi (q,\mathbf{d}) = +1.
	\end{align*}
	Then, it makes sense to talk about the \emph{$\phi$-positive} (resp. \emph{$\phi$-negative}) side of the wall $\wall$.\\
    We call $\phi$ \emph{periodic} if $\phi$ is wall-consistent and it holds:
    Let $\wall_1$ and $\wall_2$ be two parallel walls with corresponding half-spaces $\wall_1^{\varepsilon_1}$ and $\wall_2^{\varepsilon_2}$ satisfying $\wall_1^{\varepsilon_1} \subseteq \wall_2^{\varepsilon_2}$.
    Then $\wall_1^{\varepsilon_1}$ is the $\phi$-positive side of $\wall_1$ if and only if $\wall_2^{\varepsilon_2}$ is the $\phi$-positive side of $\wall_2$.
\end{defi}

One can think of several ways to introduce trivial and natural orientations on $\Sigma$.
For example, let $b$ a simplex in $\Sigma$, $\mathbf{c}$ an alcove, and $p$ a panel of $\mathbf{c}$.
Then define $\phi_b (p,\mathbf{c})=+1$ if and only if $\mathbf{c}$ and $b$ are on the same side of the wall $\wall$ containing $p$ or if $b \subseteq \wall$ to obtain the simplex orientation towards $b$, written $\phi_b$.
The reader interested in further orientations may refer to \cite[Section 3]{Graeber_2020}.
For this text, we will restrict to a specific type of (boundary) simplex orientations of affine Coxeter complexes:

\begin{defi}{(Weyl chamber orientation)} \label{def01:WeylChamberOrientation}
	Let $\Sigma$ be the Coxeter complex associated with an affine Coxeter system $(\aW,\SCox)$ and $\partial \phi_v$ the simplex orientation of the boundary complex $\partial \Sigma$ towards a chamber $v \in \partial \Sigma$.
	Then, we can assign a wall-consistent orientation to $\Sigma$ by choosing the side $\wall^\varepsilon$ of a hyperplane $\wall$ to be the positive (resp. negative) side if and only if $\partial \wall^\varepsilon$ is the positive (resp. negative) side with respect to $\partial \phi_v$.
	We denote the resulting orientation of $\Sigma$ by $\phi_v$ and call it the \emph{Weyl chamber orientation with respect to $v$}.
\end{defi}

\begin{remark}
    Using the connection between the Coxeter complex $\Sigma$ and the boundary complex $\partial \Sigma$ for affine Coxeter systems $(\aW,\SCox)$, Graeber and Schwer showed in \cite[Lemma 3.12]{Graeber_2020} that Weyl chamber orientations are wall-consistent and periodic; because of their definition via the boundary complex, they are also called \emph{induced affine orientations}.\\
    Also, observe that our definition of Weyl chamber orientations is equivalent to the following alternative formulation:
    The Weyl chamber $v \in \partial \Sigma$ corresponds to a parallel class of Weyl chambers in $\Sigma$.
    Define the orientation $\phi_v$ by setting $\phi_v (p,\mathbf{c})=+1$ if and only if $v$ has a representative cone $\cone_v$ in the same half-space defined by the supporting hyperplane $\wall_p$ of the panel $p$ as the alcove $\mathbf{c}$.
\end{remark}
%
%%%%%%%%%%%%%%%%%%%%%%%%%%%%%%%%%%%%%%%%%%%%%%%%%%%%%
%%%%%%%%%%%%%% SECTION 3 %%%%%%%%%%%%%%%%%%%%%%%%%%%%
%%%%%%%%%%%%%%%%%%%%%%%%%%%%%%%%%%%%%%%%%%%%%%%%%%%%%

\section{Positively folded galleries}\label{sec:PosFoldedGalleries}

In this section, we fix our notation on galleries in \ref{subsec:galleries}, restricting to a special case of the notion used in \cite{MilicevicSchwerThomas_2015}.
Observe that this notion is slightly different from the one used in \cite{GaussentLittelmann_2005}, where the concept of folded galleries appeared earlier. 
Further, we characterize minimal galleries by the directions of their hyperplane crossings with respect to Weyl chamber orientations (\Cref{lemma:CharacterizationOfMinimalityOfGalleriesByCrossings}) and the position of their end alcove (\Cref{lemma:CrossingDirectionsMinimalGalleries}).
In \ref{subsec:folds} below, we introduce our notation of folds and their positivity, and we also point out their connections to Coxeter shadows and Bruhat order.

\subsection{Galleries}\label{subsec:galleries}
We follow the definition of a (combinatorial) gallery given by \cite{Graeber_2020}.
To the reader familiar with (other) definitions of galleries in Coxeter systems, we want to point out that all our considered galleries are alcove-to-alcove galleries.

\begin{defi}{(Gallery)}\label{def01:AlcoveGallery}
	A \emph{gallery} in $\Sigma$ is a sequence
	\[
		\gamma = (\mathbf{c}_0, p_1, \mathbf{c}_1, \cdots, p_n, \mathbf{c}_n)
	\]
	of alcoves $\mathbf{c}_i$ and panels $p_i$ with $p_i \in \mathbf{c}_{i-1}$ and $p_i \in \mathbf{c}_i$.
    We define the \emph{length} $\length(\gamma)$ of a gallery to be $n+1$, which equals the number of alcoves contained in $\gamma$, counted with multiplicity, and we say that a gallery is \emph{minimal}, if there is no shorter gallery with start alcove $\mathbf{c}_0$ and end alcove $\mathbf{c}_n$.
\end{defi}

Often we do not require all the information on the alcoves and panels a gallery is passing, but only need to know its starting and/or ending alcoves.
For this purpose we denote abbreviatory by $\gamma: \mathbf{x} \leadsto \mathbf{y}$ a gallery $\gamma= (\mathbf{c}_0 = \mathbf{x}, p_1, \mathbf{c}_1,\dots, p_n \mathbf{c}_n = \mathbf{y})$. 
Notice that we do not call for $\mathbf{c}_i \neq \mathbf{c}_{i-1} \ \forall i \in \lbrace 1, \dots ,n \rbrace$, so that $\mathbf{c}_i = \mathbf{c}_{i-1}$ is explicitly allowed. 
Such galleries will be called \emph{stammering}, and they are not minimal.

Since a non-stammering alcove-to-alcove gallery $\gamma$ is uniquely determined by its first alcove and its panels, it makes sense to introduce the \emph{type} of a gallery as 
\[
	\tau (\gamma) := \tau (p_1) \cdots \tau (p_n) = s_{p_1} \cdots s_{p_n} \in \SCox^*.
\]
Notice that one can also introduce types for stammering galleries, but that they do not give a unique description:
Multiple stammering galleries may have the same type.
To deal with this, sometimes the notion of a \emph{decorated type} is used, marking (types of) panels between identical alcoves.
Compare \cite[Def. 4.7]{Graeber_2020}.
Since the type $\tau$ of a gallery is a word in $\SCox$ and concatenation of words gives words again, we can also define concatenation of galleries:

\begin{defi}[Concatenation of galleries]
    Let $\gamma_1= (\mathbf{c}_{0,1}, p_{1,1}, \mathbf{c}_{1,1}, \cdots, p_{n,1}, \mathbf{c}_{n,1})$ and $\gamma_2= (\mathbf{c}_{0,2}, p_{1,2}, \mathbf{c}_{1,2}, \cdots, p_{m,2}, \mathbf{c}_{m,2})$ be galleries in $\Sigma$ with $\tau(\gamma_1)= w \in \SCox^*$.
    We define the \emph{concatenation} $\tilde{\gamma} := \gamma_2 * \gamma_1$ of $\gamma_1$ and $\gamma_2$ as:
    \begin{align*}
        \tilde{\gamma}= (\mathbf{c}_{0,1}, p_{1,1}, \mathbf{c}_{1,1}, \cdots, p_{n,1}, \mathbf{c}_{n,1}, w.\mathbf{c}_{0,2}, w.p_{1,2}, w.\mathbf{c}_{1,2}, \cdots, w.p_{m,2}, w.\mathbf{c}_{m,2}).
    \end{align*}
\end{defi}

Observe that, via their type, non-stammering galleries starting in $\cf$ are in bijection to the set of words over $\SCox^*$, and minimal galleries starting in $\cf$ are in bijection to reduced expressions over $\SCox^*$.
This allows to represent, in general not uniquely, all group elements of $\aW$ by minimal galleries.
And given a Weyl chamber orientation of the Coxeter complex, we can now determine whether a given gallery is minimal by studying the crossings of the gallery.
For this, we introduce first:

\begin{defi}[Positive/negative crossing]\label{def:positiveCrossing}
    Let $(\Cox,\SCox)$ be a Coxeter system, and let the associated Coxeter complex be equipped with an orientation $\phi$.
    Consider a gallery $\gamma$ crossing a hyperplane $\wall$ at the panel $p_i$ from $\mb{c}_{i-1}$ to $\mb{c}_i$.
    We say, that the crossing of $\wall$ at $p_i$ is a \emph{positive} (resp. \emph{negative}) \emph{crossing} or \emph{in positive (negative) direction} if $\phi(p_i,\mb{c}_{i-1})=+1$ (resp. $\phi(p_i,\mb{c}_{i-1})=-1$) and $\phi(p_i,\mb{c}_{i})=-1$ (resp. $\phi(p_i,\mb{c}_{i-1})=+1$).
\end{defi}

\begin{lemma}\label{lemma:CharacterizationOfMinimalityOfGalleriesByCrossings}
    Let $(\aW,\SCox)$ be an affine Coxeter system and $x\in \aW$, $x = t^\lambda w$ with $w \in \sW$ and $ \lambda \in R^\vee$.
    Let $\phi=\phi_{v}$ be a Weyl chamber orientation of the Coxeter complex $\Sigma = \Sigma (\aW,\SCox)$ with respect to the chamber $v \in \partial \Sigma$.
    Then
    $\gamma: \cf \leadsto x\cf$ is a minimal gallery if and only if one of the following three conditions is satisfied for all positive roots $\alpha \in \Phi^+$:
    \begin{enumerate}
        \item $\langle \lambda, \alpha^\vee \rangle \neq 0$ and all crossings of hyperplanes of class $\alpha$ are either positive or all negative; 
        \item $\langle \lambda, \alpha^\vee \rangle = 0$, the alcoves $\cf$ and $\mathbf{w}$ are not separated by $\wall_{\alpha, 0}$,  and $\gamma$ crosses no $\alpha$-hyperplane;
        \item $\langle \lambda, \alpha^\vee \rangle = 0$, the alcoves $\cf$ and $\mathbf{w}$ are separated by $\wall_{\alpha, 0}$,  and $\gamma$ crosses one single $\alpha$-hyperplane once.
    \end{enumerate}
\end{lemma}

\begin{proof}
    Suppose first that $\gamma: \cf \leadsto x\cf$ is minimal.
    Then it follows from \cite[I. 4D, Prop. 4]{BrownBuildings_1989} that every hyperplane crossed by $\gamma$ is crossed only once, and the set of crossed hyperplanes is the set of hyperplanes separating $\cf$ and $\mathbf{x}$.
    We distinguish two cases:\\
    \emph{Case 1: $\langle \lambda, \alpha^\vee \rangle = 0$.}
    It follows that $x \in \mathcal{S}_{\alpha, 0}$ if $\cf$ and $\mathbf{w}$ are not separated by $\wall_{\alpha, 0}$, and that $x \in \mathcal{S}_{\alpha, -1}$, if $\cf$ and $\mathbf{w}$ are separated by $\wall_{\alpha, 0}$.
    If $x \in \mathcal{S}_{\alpha, 0}$, then $\gamma$ does not cross an $\alpha$-hyperplane, since $\cf$ is also con tained in the $0$-strip with respect to $\alpha$.
    Hence every minimal gallery $\gamma: \cf \leadsto x\cf$ satisfies $\gamma \subseteq \mathcal{S}_{\alpha,0}$ and case (2) of the assertion.
    If $x \in \mathcal{S}_{\alpha, -1}$, then $\gamma$ has to cross $\wall_{\alpha, 0}$ that separates $\cf \in \mathcal{S}_{\alpha_i, 0}$ and $x$.
    But then $\gamma$ cannot cross another $\alpha$-hyperplane without exiting the strip $\mathcal{S}_{\alpha,0}$; and with $x \in \mathcal{S}_{\alpha,0}$ it has to end in the stip and therefore had to enter it again, which contradicts the minimality of $\gamma$.
    In particular, the crossing of the hyperplane $\wall_{\alpha, 0}$ is the only crossing of an $\alpha$-hyperplane such a minimal gallery does, hence case (3) of the assertion.\\
    \emph{Case 2: $\langle \lambda, \alpha_i^\vee \rangle \neq 0$.}
    First, observe that $\cf$ and $\mb{x}$ are separated by at least one hyperplane of class $\alpha$, hence $\gamma$ has to cross $\alpha$-hyperplanes.
    Now it follows from the periodicity of $\phi$ that crossings of $\gamma$ are either all positive or all negative:
    Suppose $\gamma$ crosses not all $\alpha$-hyperplanes in the same direction, that is, one hyperplane $\wall_1$ is crossed positively and another hyperplane $\wall_2$ (with $\wall_1 = \wall_2$ not forbidden) is crossed negatively.
    Since $\mb{x}$ is on the $\varepsilon$-side and $\cf$ on the $-\varepsilon$-side of all $\alpha$-hyperplanes separating them, this implies that $\gamma$ is either crossing a hyperplane not separating $\cf$ and $\mb{x}$, or crossing one of them more than once, both items contradicting the minimality of the gallery.\\
    To prove the converse, let $\gamma$ be a gallery that satisfies one of the stated conditions for every choice of positive root $\alpha \in \Phi^+$.
    We prove via contradiction that this gallery is minimal by showing that a non-minimal gallery fails all three of the stated conditions for at least one root $\alpha_i$.
    Suppose $\gamma$ is not minimal.
    Then $\gamma$ crosses at least one hyperplane $\wall_{\alpha, k}$ with $k \in \Z$ more than once.
    Since $\phi$ is wall-consistent, a multiple crossed hyperplane is crossed at least once each in the positive and the negative direction.
    It follows immediately that $\gamma$ fails the condition of crossings of $\alpha$-hyperplanes being all positive or all negative for $\langle \lambda, \alpha^\vee \rangle \neq 0$, hence not case (1) of the assertion.
    But $\gamma$ also fails the conditions of not crossing an $\alpha$-hyperplane, hence not case (2), or only crossing a single $\alpha$-hyperplane once, hence also not case (3).
    This is a contradiction, hence $\gamma$ has to be minimal.
\end{proof}

\begin{figure}
    \centering 
    \begin{overpic}[scale=0.45]{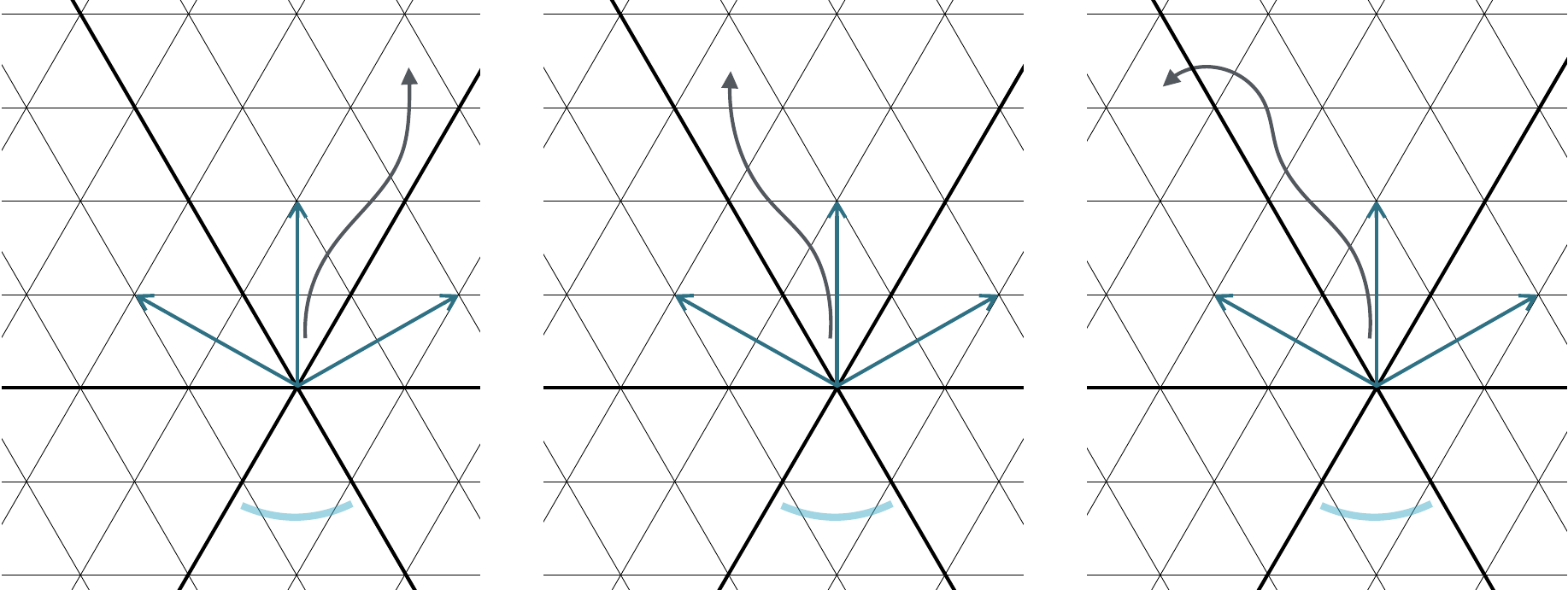}
        \put(25,15.5) {\tiny \textcolor{darkturquoise}{$\alpha_1$}}
        \put(59.5,15.5) {\tiny \textcolor{darkturquoise}{$\alpha_1$}}
        \put(94,15.5) {\tiny \textcolor{darkturquoise}{$\alpha_1$}}
    \end{overpic}
    \caption{Examples of minimal galleries satisfying the three conditions of \Cref{lemma:CharacterizationOfMinimalityOfGalleriesByCrossings} with respect to the root $\alpha_1$. The blue arc indicates the orientation defining Weyl chamber.}
    \label{fig:MinimalGalleriesWithDifferentCrossings}
\end{figure}

\begin{example}
    \Cref{fig:MinimalGalleriesWithDifferentCrossings} shows examples of the three possible situations of crossings discussed in \Cref{lemma:CharacterizationOfMinimalityOfGalleriesByCrossings}.
    The picture on the left shows a gallery ending in $x=t^\lambda w$ with $\langle \lambda, \alpha_1^\vee \rangle = 3$.
    It crosses two $\alpha_1$-hyperplanes, with the second and the fourth panels.
    Both crossings are in the same direction.\\
    The picture in the middle shows a minimal gallery ending in $x=t^\lambda w$ with $\langle \lambda, \alpha_1^\vee \rangle = 0$ and $\mb{w} \subseteq \wall_{\alpha_1, 0}^{\cf}$.
    That is $\gamma \subseteq \mathcal{S}_{\alpha_1,0}$, hence $\gamma$ crosses no $\alpha_1$-hyperplane.\\
    The picture on the right shows a minimal gallery ending in $x= t^\lambda w$ with $\langle \lambda, \alpha_1^\vee \rangle = 0$ and $\mb{w} \nsubseteq \wall_{\alpha_1, 0}^{\cf}$.
    The gallery crosses $\wall_{\alpha_1,0}$ and no other $\alpha_1$-hyperplane.
\end{example}

We now want to study the directions of crossings with respect to a given Weyl chamber orientation in more detail.
Given $\lambda \in R^\vee$ and $\alpha \in \Phi^+$, observe that $\sgn(\langle \lambda, \alpha^\vee\rangle)$ only depends on the spherical Weyl chamber $\mathcal{C}_w$ which contains the lattice point $\lambda$.
Hence the direction of crossing $\alpha$-hyperplanes with respect to a Weyl chamber orientation is the same for all minimal galleries in this chamber and can be deduced by checking if $\mathbf{w}$ is on the same side of $\mathcal{H}_{\alpha,0}$ as $\cf$.
We will discuss this in the next lemma:

\begin{lemma}[Crossing directions of minimal galleries] \label{lemma:CrossingDirectionsMinimalGalleries}
    Let $(\aW,\SCox)$ be an affine Coxeter system with corresponding spherical group $\sW$.
    Denote by $\Sigma = \Sigma (\aW,\SCox)$ its Coxeter complex and let $\phi=\phi_{w}$ be the Weyl chamber orientation of $\Sigma$ with respect to $w \in \sW$.
    Let further $\partial \phi_{w}$ be the alcove orientation of the spherical Coxeter complex $\partial \Sigma$ that induces $\phi_{w}$.
    Then one has for all $v \in \sW$, $x \in \cone_v$, and every minimal gallery $\gamma : \cf \leadsto x\cf$ that crosses a hyperplane $\wall$:
    $\gamma$ crosses $\wall$ in $\Sigma$ positively with respect to $\phi_{w}$ if and only if
    $v$ is on the $\partial \phi_{w}$-negative side of the hyperplane $\partial \wall$ in $\partial \Sigma$.
\end{lemma}

\begin{proof}
    Recall that since $\gamma$ is minimal, no hyperplane is crossed more than once by the gallery and it is reasonable to talk about \emph{the} crossing of the hyperplane.\\
    First, let $\wall_{\alpha, k}$ with $\alpha \in \Phi^+$, $k \in \Z$ be a hyperplane $\gamma$ crosses in $\phi_w$-positive direction in $\Sigma$.
    It follows that $x \in \wall_{\alpha, k}^-$, and with \Cref{lemma:CharacterizationOfMinimalityOfGalleriesByCrossings}, condition (1) or (3), also $x \in \wall_{\alpha, 0}^-$.
    Since $\cone_v$ is bounded by $\wall_{\beta, 0}$ for suitable $\beta \in \Phi^+$, it holds either $\cone_v \subseteq \wall_{\alpha,0}^+$ or $\cone_v \subseteq \wall_{\alpha,0}^-$.
    But we have seen that $x \in \wall_{\alpha, 0}^-$, and with $x \in \cone_v$ it follows $\cone_v \subseteq \wall_{\alpha, 0}^-$.
    With $\mathbf{v} \subseteq \cone_v$, the claim now follows from the definition of induced spherical orientation on $\partial \Sigma$.\\
    On the other hand, if $v$ is on the $\partial \phi_w$-negative side of a hyperplane $\partial \wall$ in $\partial \Sigma$, it follows from the definition of induced affine orientation that $\cone_v \subseteq \wall_{\alpha, 0}^-$, with $\alpha$ the parallelism class of the hyperplanes in $\Sigma$ that belong to the hyperplane $\partial \wall$ in the boundary.
    Since $x \in \cone_v$, it holds $x \in \wall_{\alpha, 0}^-$.
    And since $\gamma$ is minimal, with \Cref{lemma:CharacterizationOfMinimalityOfGalleriesByCrossings}, case (1) or (3), we further have $x \in \wall_{\alpha, k}^-$ for all $\alpha$-hyperplanes in $\Sigma$ that $\gamma$ crosses.
\end{proof}

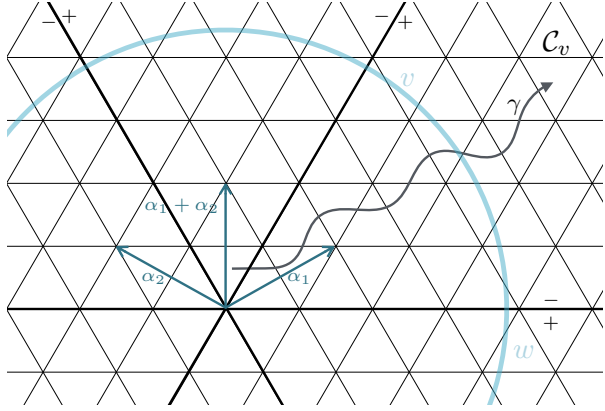
\begin{figure}
    \centering 
    \begin{overpic}[scale=0.45]{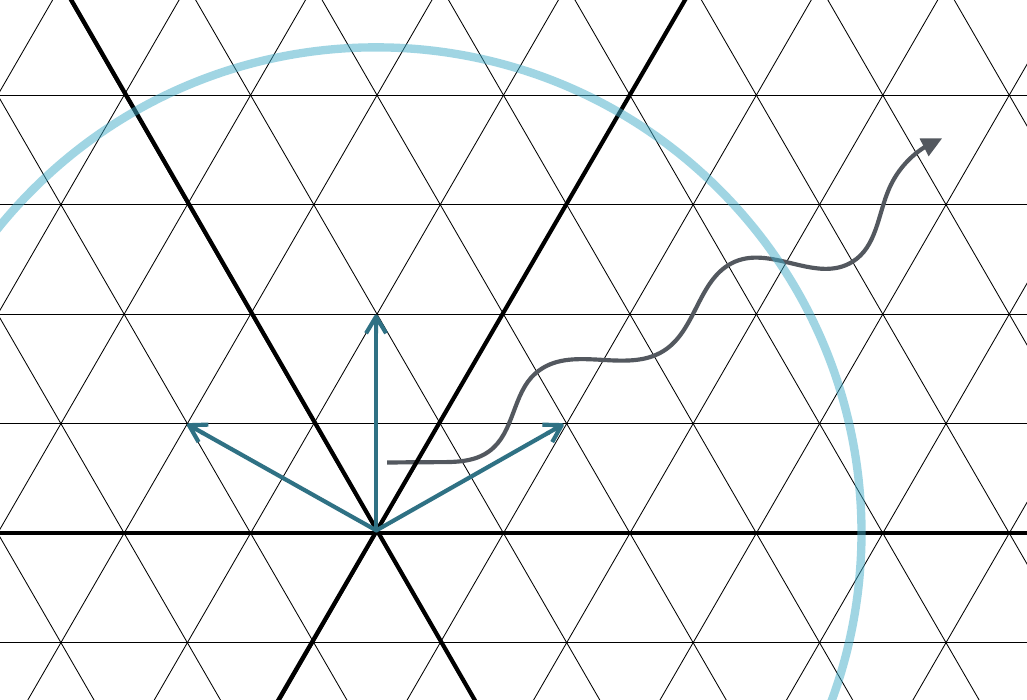}
        \put(22.5,21) {\tiny \textcolor{darkturquoise}{$\alpha_2$}}
        \put(47,21) {\tiny \textcolor{darkturquoise}{$\alpha_1$}}
        \put(23,33) {\tiny \textcolor{darkturquoise}{$\alpha_1+\alpha_2$}}
        \put(84,50) {\footnotesize $\gamma$}
        \put(90,60) {$\cone_v$}
        \put(66,54) {\textcolor{pastelblue}{$v$}}
        \put(85,8) {\textcolor{pastelblue}{$w$}}
        %orientation
        \put(90,13) {\tiny $+$}
        \put(90,17) {\tiny $-$}
        \put(65.5,64) {\tiny $+$}
        \put(61.5,65) {\tiny $-$}
        \put(9,65) {\tiny $+$}
        \put(5.5,64) {\tiny $-$}
    \end{overpic}
    \caption{Visualization of \Cref{lemma:CrossingDirectionsMinimalGalleries}. The light blue circle indicates the boundary sphere $\partial \Sigma$. The depicted gallery $\gamma$ ends in $\cone_v$. Since $v \in \partial \Sigma$ is on the $\partial \phi_w$-negative side of the horizontal hyperplane $\partial \wall$, crossings of horizontal hyperplanes of $\gamma$ in $\Sigma$ are positive. Since $\partial v$ is on the $\partial \phi_w$-positive side of the boundary hyperplane corresponding to hyperplanes in $\Sigma$ perpendicular to $\alpha_1$, all crossings of $\alpha_1$-hyperplanes of $\gamma$ are $\phi_w$-negative.}
    \label{fig:CrossingDirectionsMinimalGalleries}
\end{figure}

See \Cref{fig:CrossingDirectionsMinimalGalleries} for an example.
The result of \Cref{lemma:CrossingDirectionsMinimalGalleries} also applies to minimal galleries not starting in the fundamental alcove:

\begin{cor}\label{cor:CrossingDirectionsMinimalGalleries--NoncanonicGalleries}
    Let $(\aW,\SCox)$ be an affine Coxeter system with corresponding spherical group $\sW$.
    Denote by $\Sigma = \Sigma (\aW,\SCox)$ its Coxeter complex and let $\phi=\phi_{w}$ be the Weyl chamber orientation of $\Sigma$ with $w \in \sW$.
    Let further $\partial \phi_{w}$ be the alcove orientation of the spherical Coxeter complex $\partial \Sigma$ that induces $\phi_{w}$.
    Then one has for all $u,v \in \sW$, $y=t^\mu u \in \aW$, $x \in \cone_{\mu,v}$ and every minimal gallery $\gamma : y\cf \leadsto x\cf$ that crosses a hyperplane $\wall$:
    The crossing is $\phi_w$-positive if and only if
    $v$ is on the $\partial \phi_{w}$-negative side of the hyperplane $\partial \wall$ in $\partial \Sigma$.
\end{cor}

\begin{proof}
    Let $\gamma: y\cf \leadsto x\cf$ be a minimal gallery of type $s_{\gamma_1}s_{\gamma_2}\cdots s_{\gamma_n}$ and $\setwalls_\gamma$ the set of crossed hyperplanes.
    Assume that $y=t^0 u$, that is, one vertex of the alcove $\mathbf{y}=y\cf$ is the origin.
    Then $x \in \cone_{0,v}=\cone_v$.
    Let $\gamma_y : \cf \leadsto \mathbf{y}$ be a minimal gallery with corresponding word $y= s_{y_1}s_{y_2}\cdots s_{y_m}$ and set of crossed hyperplanes $\setwalls_y$.
    We distinguish two cases:\\
    \emph{Case 1: $\setwalls_y \cap \setwalls_\gamma = \emptyset$.}
    Then the concatenation $\gamma' := \gamma * \gamma_y: \cf \leadsto x\cf$ is a minimal gallery.
    Application of \Cref{lemma:CrossingDirectionsMinimalGalleries} proves the assertion.\\
    \emph{Case 2: $\setwalls_y \cap \setwalls_\gamma \neq \emptyset$.}
    Then the concatenation $\gamma' := \gamma * \gamma_y: \cf \leadsto x\cf$ is not a minimal gallery.
    We observe:
    For a hyperplane $\wall$ the following two assertions hold:
    \begin{itemize}
        \item $\wall \in \setwalls_y$ if and only if $\cf$ and $\mathbf{y}$ are on different sides of $\wall$ and $\wall$ contains the origin;
        \item $\wall \in \setwalls_\gamma$ if $\mathbf{y}$ and $\mathbf{v}$ are on different sides of $\wall$; with equivalence if $\wall$ contains the origin.
    \end{itemize}
    Hence for $\wall \in \setwalls_y \cap \setwalls_\gamma$, we have that $\cf$ and $\mathbf{v}$ are on one side and $\mathbf{y}$ is on the other side of $\wall$.
    And for hyperplanes $\wall$ with $\wall \in \setwalls_\gamma$ and $\wall \notin \setwalls_y$ containing the origin, we have that $\cf$ and $\mathbf{y}$ are together on one side of $\wall$ and $\mathbf{v}$ is on the other side of $\wall$.
    Given this, let now $\wall$ be a hyperplane crossed by $\gamma$. 
    We consider the two possible situations:\\
    \emph{Case 2A: $\wall \in \setwalls_y \cap \setwalls_\gamma$.}
    Then $\wall$ is crossed twice by $\gamma'$, once in every concatenated part.
    With our observations from above and the periodicity of the orientation $\phi_w$, this gives us that $\gamma$ is crossing $\wall$ $\phi_{w}$-positively if and only if $\gamma_y$ crosses $\wall$ $\phi_{w}$-negatively.
    With \Cref{lemma:CrossingDirectionsMinimalGalleries} this is equivalent to $\zeta(y)=u$ being on the $\partial \phi_w$-positive side of $\partial \wall$ in $\partial \Sigma$.
    And again with our observations above this holds if and only if $v$ is on the $\partial \phi_w$ negative side of $\partial \wall$ in $\partial \Sigma$, hence the claim.\\
    \emph{Case 2B: $\wall \in \setwalls_\gamma$ and $\wall \notin \setwalls_y$.}
    It follows from our observation above on $\cf$ and $\mathbf{y}$ being on the same side of $\wall$, that a minimal gallery $\tilde{\gamma}: \cf \leadsto x\cf$ crosses $\wall$ in the same direction as $\gamma$.
    Hence  $\gamma$ crosses $\wall$ $\phi_{w}$-positively if and only if $\tilde{\gamma}$ crosses $\wall$ $\phi_{w}$-negatively.
    Application of \Cref{lemma:CrossingDirectionsMinimalGalleries} to $\tilde{\gamma}$ gives:
    $\tilde{\gamma}$ crosses $\wall$ $\phi_{w}$-positively if and only if $v$ is on the $\partial \phi_w$ negative side of $\partial \wall$ in $\partial \Sigma$.
    This proves the assertion.\\
    To see now that this also holds for $y=t^\mu u$ with $\mu \neq 0$, let $\gamma_t: \mathbf{y} \leadsto t^0u$ be a minimal gallery and consider $\gamma'' := \gamma * \gamma_t = t^{-\mu}(\gamma)$.
    Since $\phi_w$ is periodic, $\gamma''$ crosses all hyperplanes in the $\gamma$-part in the same direction as $\gamma$.
    Application of the results above to this part proves the claim.
\end{proof}

From \Cref{lemma:CrossingDirectionsMinimalGalleries} it follows, that all minimal galleries with end alcoves in a common Weyl chamber show the same behavior in terms of crossing directions of parallel hyperplanes with respect to a Weyl chamber orientation $\phi_w$.
This will be crucial when it comes to positively folding these galleries as we will see in the next subsection and \Cref{thm:positiveFoldingPatternsAndMomentGraph} below.

\subsection{Folds}\label{subsec:folds}

We now introduce folded galleries:

\begin{defi}[Folded gallery]
    Let $\gamma= ( \mathbf{c}_0, p_1, \mathbf{c}_1, \cdots, p_n \mathbf{c}_n)$ be a gallery.
    If $\mathbf{c}_i = \mathbf{c}_{i-1}$ for some $i$, we say that a gallery is \emph{stammering} and \emph{folded at the panel} $p_i$, otherwise, we call it \emph{unfolded}.
\end{defi}

When only considering the end alcove, folding a gallery $\gamma$ at a panel $p_i$ is equivalent to deleting the letter $\tau(p_i)= s_{p_i}$ in the word for $\gamma$.
But notice that every folded gallery $\gamma'$, obtained from an unfolded gallery $\gamma$, satisfies: $\tau(\gamma')= \tau(\gamma)$.
Geometrically, folding a gallery $\gamma$ at a panel $p_i$ reflects the remaining part of the gallery $\gamma_r := (\mathbf{c}_i, p_{i+1}, \cdots, p_n, \mathbf{c}_n)$ at the hyperplane $\wall$ containing the panel $p_i$.
The following definition of this explicit (un-)folding of a gallery is taken from \cite[Def. 4.15]{Graeber_2020} and states this more precisely:

\begin{defi}{((Un-)Folding a gallery)}\label{def:(un)foldingAGallery}
	Let $\gamma = ( \mathbf{c}_0, p_1, \mathbf{c}_1, \dots, p_n, \mathbf{c}_n)$ be a gallery and denote by $r_i$ the reflection across the wall $\mathcal{H}_{p_i}$ spanned by the panel $p_i$.
	Then we call
	\begin{align*}
		\gamma^i := (\mathbf{c}_0, p_1, \mathbf{c}_1, \dots, p_i, r_i\mathbf{c}_i, r_ip_{i+1}, r_i\mathbf{c}_{i+1}, \dots, r_i p_n, r_i \mathbf{c}_n)
	\end{align*}
	a \emph{(un-)folding of $\gamma$ at the panel $p_i$}.
\end{defi}

We want to study folded galleries in relation to orientations of $\Sigma$:

\begin{defi}[Positive/negative folds]
	Let $\phi$ be an orientation of the Coxeter complex $\Sigma$ and $\gamma = (\mathbf{c}_0, p_1, \dots, p_n, \mathbf{c}_n)$ be a gallery in $\Sigma$.
    We say, that $\gamma$ is \emph{positively} (resp. \emph{negatively}) \emph{folded with respect to the orientation $\phi$ at the panel $p_i$} if $\mathbf{c}_i = \mathbf{c}_{i-1}$ and $\phi (p_i, \mathbf{c}_i)= +1$ (resp. $\phi ( p_i, \mathbf{c}_i)= -1$).
	We say further, that $\gamma$ is \emph{positively} (resp. \emph{negatively}) \emph{folded with respect to the orientation $\phi$} if for all $i \in \lbrace 1, \dots, n \rbrace$ either $\mathbf{c}_i \neq \mathbf{c}_{i-1}$, that is $\gamma$ is unfolded, or $\mathbf{c}_i = \mathbf{c}_{i-1}$ and $\phi (p_i, \mathbf{c}_i)= +1$ (resp. $\phi ( p_i, \mathbf{c}_i)= -1$).
\end{defi}

Note that it is possible to fold a gallery at several panels simultaneously and that foldings commute: $\left( \gamma^i\right)^j = \left(\gamma^j\right)^i$ and $\left(\gamma^i\right)^i = \gamma$.
This was proven in \cite[Lemma 4.17]{Graeber_2020}.
Denote the set of indices of panels, where $\gamma$ has a fold in, by $F(\gamma)$.
Since reflections are type preserving on $\Sigma$, it holds $\tau (\gamma) = \tau (\gamma^i)$.
Depending on whether $\gamma$ was already folded at the panel $p_i$, the \emph{number of folds} $\labs F(\gamma)\rabs$ of the gallery $\gamma$ increases or decreases by one for each folding.

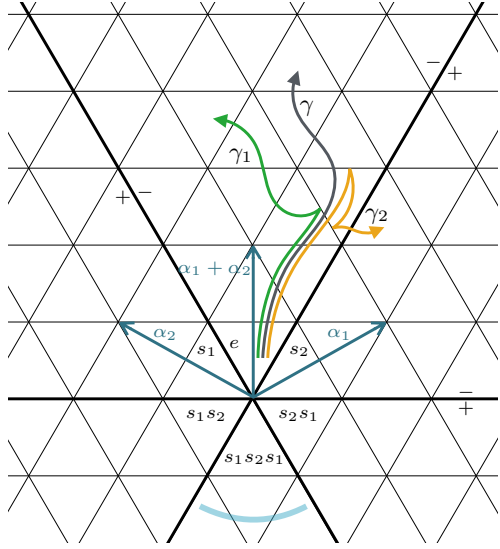
\begin{figure}
    \centering 
    \begin{overpic}[scale=0.55]{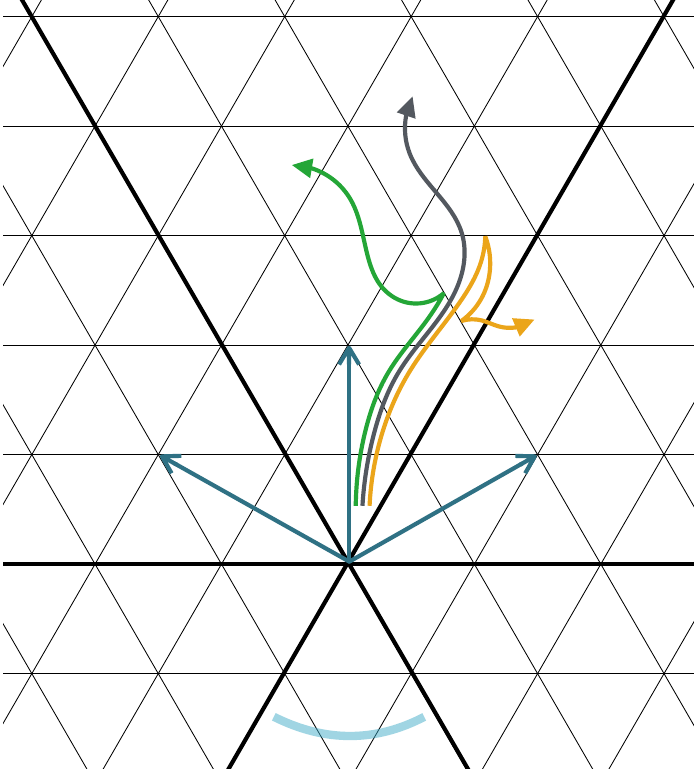}
        \put(27,38) {\tiny \textcolor{darkturquoise}{$\alpha_2$}}
        \put(59,38) {\tiny \textcolor{darkturquoise}{$\alpha_1$}}
        \put(31.5,50) {\tiny \textcolor{darkturquoise}{$\alpha_1+\alpha_2$}}

        \put(41,36) {\tiny $e$}
        \put(35,35) {\tiny $s_1$}
        \put(33,23) {\tiny $s_1s_2$}
        \put(40,15) {\tiny $s_1s_2s_1$}
        \put(50,23) {\tiny $s_2s_1$}
        \put(52,35) {\tiny $s_2$}

        \put(54,80) {\footnotesize $\gamma$}
        \put(41,71) {\footnotesize $\gamma_1$}
        \put(66,60) {\footnotesize $\gamma_2$}
         %orientation
        \put(83,24) {\tiny $+$}
        \put(83,27) {\tiny $-$}
        \put(81,86) {\tiny $+$}
        \put(77,88) {\tiny $-$}
        \put(20,63) {\tiny $+$}
        \put(24,64) {\tiny $-$}
     \end{overpic}
    \caption{Folded galleries.}
    \label{fig:FoldedGalleries}
\end{figure}

\begin{example}
    \Cref{fig:FoldedGalleries} shows two galleries $\gamma_1$ and $\gamma_2$ obtained by folding the grey gallery $\gamma$.
    $\gamma_1$ is folded once at the fourth panel.
    With respect to the Weyl chamber orientation defined by $s_1s_2s_1=w_0 \in \partial \Sigma$, this fold is positive, since the repeated alcove is on the positive side of the $\alpha_1$-hyperplane containing the fourth panel.
    $\gamma_2$ is folded twice at the fifth and sixth panels.
    The fold at the fifth panel is positive with respect to $\phi_{w_0}$, but the second fold is negative:
    The repeated sixth and seventh alcove of $\gamma_2$ is on the negative side of the $\alpha_1$-hyperplane containing the sixth panel.
    Hence $\gamma_2$ is not positively folded.
\end{example}

The positivity of a fold with respect to a Weyl chamber orientation $\phi_w$ only depends on the parallelism class of the hyperplane the folds panel is contained in.
With our results from \Cref{lemma:CrossingDirectionsMinimalGalleries} above, that galleries ending in the same Weyl chamber cross parallel hyperplanes in the same crossing direction, we conclude that the set of hyperplane classes these galleries can be folded in positively must coincide.
Hence we introduce the following:

\begin{defi}[(Positive) Folding pattern]
    Let $\gamma$ be a minimal gallery with $n$ folds. If the $i$-th fold, counted with a rising index number, is in a panel that has supporting hyperplane in the parallelism class $\alpha_j$ for a positive root $\alpha_j$, then we call the $n$-tuple $(\alpha_{j,i})_{i = \lbrack n \rbrack}$ the \emph{folding pattern} of the gallery $\gamma$.
    If the Coxeter complex $\Sigma=\Sigma (\Cox,\SCox)$ is equipped with a periodic orientation $\phi$ and all foldings in a folding pattern  $(\alpha_{j,i})_{i = \lbrack n \rbrack}$ are positive with respect to $\phi$, we call this a \emph{$\phi$-positive folding pattern}.
\end{defi}

\begin{remark}
    Restricting to periodic orientations of $\Sigma$ when defining positive folding patterns allows the application of a pattern to arbitrary crossed hyperplanes of the given parallelism class for a fixed gallery $\gamma$ while keeping the positivity.
\end{remark}
%
%%%%%%%%%%%%%%%%%%%%%%%%%%%%%%%%%%%%%%%%%%%%%%%%%%%%%%%%%%
%%%%%%%%%%%%%%%%%%%%%%%%%%%%%%%%%%%%%%%%%%%%%%%%%%%%%%%%%%
%%%%%% SECTION 4 %%%%%%%%%%%%%%%%%%%%%%%%%%%%%%%%%%%%%%%
%%%%%%%%%%%%%%%%%%%%%%%%%%%%%%%%%%%%%%%%%%%%%%%%%%%%%%%%%%

\section{Connection to moment graphs}\label{sec:momentGraphs}

In this section, we point out two connections between the Bruhat moment graph of the spherical Weyl group $\sW$ associated with an affine Coxeter group $\aW$ and the $\phi$-positive folding patterns in $\aW$ for a Weyl chamber orientation $\phi_{w}$.
We start by giving a general definition of a moment graph, following \cite{Fiebig2016}:

\begin{defi}[Moment graph]
    Let $\Lambda = \Z^r$ be a lattice.
    A moment graph $\mathcal{G}$ over $\Lambda$ is a labeled simple directed graph $\graph= (V,E)$ with vertices $V$, directed edges $E$, and labeling map $f_\graph: E \to \Lambda \setminus \lbrace 0 \rbrace$.
\end{defi}

The class of moment graphs important to us are moment graphs obtained from Bruhat graphs.
Bruhat graphs visualise Bruhat order on Coxeter groups.
It is well known that Bruhat order has a geometric interpretation via galleries: For group elements $x,y \in \Cox$ it holds $y \leq x$ in Bruhat order if and only if a minimal gallery $\gamma_x: \cf \leadsto x\cf$ can be folded onto a gallery $\gamma_y: \cf \leadsto y\cf$ of the same type.
See \cite[Prop. 6.8]{Graeber_2020}.
We will prove another connection between Bruhat graphs and positive folding patterns of galleries below.
For this, let us recall the definition of a Bruhat graph from \cite[Section 2]{BjornerBrenti_2006}:

\begin{defi}[Bruhat graph]\label{def:BruhatGraph}
    Let $(\Cox,\SCox)$ be a Coxeter system and $R = \lbrace w s w^{-1}:\ w \in \Cox, \ s\in \SCox \rbrace$ its set of reflections.
    The \emph{Bruhat graph} of $\Cox$ is the directed graph $\mathcal{B}_\Cox= (V, E)$ with $V= \lbrace w:\ w \in \Cox \rbrace$ and $e = (u,w) \in E \Leftrightarrow  l(u)< l(w)$ and $w=r  (u) = u \cdot r$ for a reflection $r \in R$.
\end{defi}

Since there is a 1-on-1 correspondence between reflections in $\sW$ and the positive roots of the affine Coxeter group $\aW$, each reflection $r$ corresponds to one parallelism class of hyperplanes in the affine Coxeter group.
Therefore we can index the reflections with their affiliated positive root in $\aW$ and obtain a moment graph:

\begin{defi}[Bruhat moment graph]
    Let $(\sW,\SCox_0)$ be a spherical Coxeter system with root system $\Phi$ and $\Phi^+$ the set of positive roots.
    Denote by $X:= \lbrace \lambda \in T \ \vert \ \langle \lambda, \alpha^\vee \rangle \in \Z \ \forall \alpha \in \Phi \rbrace$ the weight lattice.
    The labeled graph $\graph_{\sW}= (V,E, \Phi^+)$ with vertex set $V= \lbrace w:\ w \in \sW \rbrace$ and edges $e = (u,w) \in E$ if and only if $l(u)< l(w)$ and $w=r_{\alpha_i} (u)$ for a reflection $r_{\alpha_i} \in R$ with $\alpha_i \in \Phi^+$ a positive root in $\sW$ is called the \emph{Bruhat moment graph}.
    We set $f_{\graph_{\sW}} (e= (u,w))= \alpha_i$.
\end{defi}

For example, the picture on the right in \Cref{fig:BruhatMomentGraphAndPosFolding} depicts the Bruhat moment graph of the Coxeter group of type $\tA_2$.

\subsection{Folding patterns and moment graphs}

Having this, we are ready to state our main result:
The Bruhat graph encodes adjacency relations via reflections of group elements in $\sW$, and the labeling of the Bruhat moment graph gives the parallelism classes of the hyperplanes corresponding to these reflections.
Because of the semidirect product structure of affine Coxeter groups, the reflections in $\sW$ encoded in the Bruhat graph each represent a class of reflections in parallel hyperplanes in the associated affine group $\aW$.
For periodic orientations of its Coxeter complex $\Sigma$, it then suffices to study a single representative of a hyperplane class to understand the positivity of foldings in affine Coxeter groups, hence a Bruhat moment graph.
This gives the following connection:

\begin{theorem}[Positive folding patterns and moment graph]\label{thm:positiveFoldingPatternsAndMomentGraph}
    Let $(\aW,\SCox)$ be an affine Coxeter system and $\sW$ the corresponding spherical Coxeter group with moment graph $\graph_{\sW}= (V, E, \Phi^+_{\sW})$.
    Denote by $\Sigma$ the affine Coxeter complex corresponding to $(\aW,\SCox)$ and let $\phi_{w_0}$ be the Weyl chamber orientation on $\Sigma$ with $w_0 \in \sW$ the unique longest element.
    Let $x \in \cone_v$ and $\gamma: \cf \leadsto x\cf$ be a minimal gallery.
    Then $(\alpha_{j,i})_{i = \lbrack n \rbrack}$ 
   is a positive folding pattern for $\gamma$ if and only if $(\alpha_{j,1}, \alpha_{j,2}, \dots, \alpha_{j,n})$ is the label index sequence of a directed path in $\graph_{\sW}$ starting in the vertex $v$.  
\end{theorem}

Before we present the proof, we recover a well-known result on the maximum number of positive folds of a minimal gallery with respect to a Weyl chamber orientation from this:

\begin{cor}[\cite{Graeber_2020}]
    Let $(\aW,\SCox)$ be an affine Coxeter system and $\sW$ the corresponding spherical Coxeter group with unique longest element $w_0 \in \sW$.
    Let further $\gamma$ be a minimal gallery that is positively folded with respect to the Weyl chamber orientation $\phi_{w_0}$.
    Then one has for the number of folds $\labs F(\gamma)\rabs$:
    \begin{align*}
        \labs F(\gamma) \rabs \leq l(w_0),
    \end{align*}
    where $l$ denotes the word length of an element in $\sW$.
\end{cor}

\begin{proof}
    This follows immediately from \Cref{thm:positiveFoldingPatternsAndMomentGraph} since no directed path in the Bruhat moment graph is longer than $l(w_0)$.
\end{proof}

Next, we want to describe a connection between directed edges in the moment graph and the Weyl chamber orientation $\phi_{w_0}$ for $w_0 \in W_0$ as a preparation for the proof of \Cref{thm:positiveFoldingPatternsAndMomentGraph}:

\begin{lemma}[Bruhat moment graph and $\phi_{w_0}$]\label{lemma:EdgesInMomentGraphAndW0Orientation}
    Let $\graph_{\sW}$ be the moment graph corresponding to the spherical Coxeter group $\sW$ and affine Coxeter system $(\aW,\SCox)$ and let $u,w \in \sW$.
    Denote by $\partial \Sigma$ the Coxeter complex associated with $\sW$ in the boundary of $\Sigma$ and by $\partial\phi_{w_0}$ the Weyl chamber orientation of $\partial\Sigma$ induced by the unique longest element $w_0 \in \sW$.
    Then there exists a directed edge $e= (u,w) \in \graph_{\sW}$ if and only if the Weyl chamber $u \subseteq \partial \Sigma$ is on the $\partial\phi_{w_0}$-negative side and $w \subseteq \partial \Sigma$ is on the $\partial\phi_{w_0}$-positive side of the reflection hyperplane $\wall_r$ with $w= r(u) =  u \cdot r$, $r \in R$.
\end{lemma}

\begin{proof}
    First, let $e=(u,w) \in \graph_{\sW}$ be a directed edge from $u$ to $w$ in the moment graph of $\sW$.
    It follows that $w= u \cdot r$ for a reflection $r \in R$ and $l(u) < l(w)$.
    Observe with \cite[Section 1.8]{Humphreys_1990}, that
    \begin{align*}
        l(u)< l(w) &\Leftrightarrow -l(u) > -l(w) \Leftrightarrow l(w_0) - l(u) > l(w_0)- l(w) \Leftrightarrow l(u^{-1}w_0)> l(w^{-1}w_0) \\
        &\Leftrightarrow \d(u,w_0) > \d(w,w_0).
    \end{align*}
    Let $\partial \wall_r$ be the hyperplane in the Coxeter complex $\partial \Sigma$ corresponding to $r$.
    $u$ and $w$ lie on opposite sides of $\partial \wall_r$ and it follows from the definition of $\phi_{w_0}$ that $w_0$ is on the $\phi_{w_0}$-positive side of $\partial \wall_r$.
    Hence if $w=w_0$, $w$ is on the $\phi_{w_0}$-positive side of $\partial \wall_r$ and there is nothing more to show.
    If $w \neq w_0$, we show that $w$ is on the $\phi_{w_0}$-positive side of $\partial \wall_r$ via contradiction.
    Suppose that $w$ is on the $\phi_{w_0}$-negative side of $\partial \mathcal{H}_r$.
    Then $u$ is on the $\phi_{w_0}$-positive side of $\partial \wall_r$.
    And since $w_0$ is on the $\phi_{w_0}$-positive side of all hyperplanes in $\partial \Sigma$, this implies $\d (u, w_0) < \d(w,w_0)$, because:
     Let $\gamma: w_0 \leadsto w$ be a minimal gallery in $\partial \Sigma$.
     Then we have $\length (\gamma) = \d(w, w_0)$.
     Since $\partial \wall_r$ separates $w_0$ and $w$, $\gamma$ crosses $\partial \wall_r$ with a panel $p$.
     Now fold $\gamma$ at $p$ to obtain the gallery $\gamma'$.
     Then $\gamma'$ is stuttering, hence not minimal, fully contained in the half-space defined by $\partial \wall_r$ that contains $w_0$ and, since $w = u \cdot r$, $\gamma' : w_0 \leadsto u$.
     But then we can shorten $\gamma'$ to obtain a minimal gallery $\gamma'': w_0 \leadsto u$ with $\length (\gamma'') < \length(\gamma')= \length (\gamma)$.
     With $\d(u,w_0)= \length (\gamma'')$ follows the claim: 
     $\d(u,w_0) < \d(w, w_0)$.
    But this is a contradiction.
    Hence $w$ is on the $\phi_{w_0}$-positive and $u$ on the $\phi_{w_0}$-negative side of $\partial \wall_r$.\\
    To prove the converse, let $u$ be on the $\phi_{w_0}$-negative and $w$ be on the $\phi_{w_0}$-positive side of the hyperplane $\partial \wall_r$ with $w = u \cdot r$ and $r \in R$.
    Then surely $(u,w)$ or $(w,u)$ is an edge of the moment graph $\graph_{\sW}$, depending on whether $l(u)< l(w)$ or $l(w) < l(u)$.
    Since $w$ is on the $\phi_{w_0}$-positive side of the hyperplane $\partial \wall_r$, it follows from the definition of the orientation $\partial \phi_{w_0}$ that $w$ and $w_0$ lie on the same side of the hyperplane $\partial \wall_r$.
    Because $\partial \wall_r$ separates $w$ and $u$, we see with the same argumentation as above that $\d(w,w_0) < \d(u,w_0)$.
    And with this, we have with \cite[Section 1.8]{Humphreys_1990} again:
    \begin{align*}
        \d(w,w_0) < \d(u,w_0) 
        &\Leftrightarrow l(u) < l(w).
    \end{align*}
    Hence $(u,w) \in \graph_{\sW}$.
\end{proof}

Having this, we are ready to prove \Cref{thm:positiveFoldingPatternsAndMomentGraph}:

\begin{proof}[Proof of \Cref{thm:positiveFoldingPatternsAndMomentGraph}.]
    First, let $(\alpha_{j,i})_{i = \lbrack n \rbrack}$ be a positive folding pattern for $\gamma$.
    Suppose, $(\alpha_{j,i})_{i = \lbrack n \rbrack}$ is not a label index sequence of a directed path in $\graph_{\sW}$ starting in $v$.
    Then we find an index $j,k$ with $k \in \lbrace 1, \dots, n \rbrace$, such that $(\alpha_{j,1}, \dots, \alpha_{j,{k-1}})$ is an index sequence of a directed path $\xi = (v, \dots, u)$ in $\graph_{\sW}$ with $v, u$ vertices, and $(\alpha_{j,1}, \dots, \alpha_{j,k})$ is not.
    Consider the gallery $\gamma'$ which is obtained from $\gamma$ by application of the folding pattern $(\alpha_{j,i})_{i = \lbrack k-1 \rbrack}$ to a feasible choice of panels, such that $\gamma'$ crosses at least one hyperplane of class $\alpha_{j,k}$ in the remaining unfolded part behind the last fold.
    If we can not find such a gallery, the folding pattern $(\alpha_{j,i})_{i = \lbrack k \rbrack}$ does not apply to $\gamma$, and since the pattern $(\alpha_{j,i})_{i = \lbrack k-1 \rbrack}$ corresponded to a directed path in the graph $\graph_{\sW}$ starting in $v$ there is nothing to show.
    So suppose we find such a gallery $\gamma'$.
    Denote the unfolded remaining part of the gallery $\gamma'$ by $\gamma'_r: y\cf \leadsto z \cf$ with $y = t^\mu w$ and $z= t^\nu w'$.
    Since $\gamma'_r$ is minimal and assumed to cross at least one hyperplane $\wall_{\alpha_{i,k}, h}$ with $\alpha_{i,k}$ the next element in the positive folding pattern for $\gamma$, we know that $z\cf$ satisfies $z\cf \subseteq \wall_{\alpha_{i,k}, h}^-$ with respect to the orientation $\phi_{w_0}$.
    Hence we conclude with \Cref{cor:CrossingDirectionsMinimalGalleries--NoncanonicGalleries} that $z\cf \subseteq \mathcal{C}_{\mu, u'}$ with $u'$ on the $\partial \phi_{w_0}$-negative side of the hyperplane $\partial \wall_{\alpha_{i,k},h}$ in $\partial \Sigma$.
    But then \Cref{lemma:EdgesInMomentGraphAndW0Orientation} shows, that there is a directed edge $e$ labeled $\alpha_{i,k}$ leaving the node $u$ in the moment graph.
    So we can prolong the directed path $\xi$ above with $e$ such that the label index sequence $(\alpha_{i,1}, \dots, \alpha_{i,{k}})$ does correspond to a directed path in the moment graph $\graph_{\sW}$.\\
    To prove the converse, let $(\alpha_{v_1}, \alpha_{v_2}, \dots, \alpha_{v_n})$ be the label index sequence of a directed path in $\graph_{\sW}$ starting in the node labeled $v$.
    Then $v$ is on the $\partial \phi_{w_0}$-negative side of the hyperplanes $\partial \wall_{\alpha_{v_i}}$ for $i \in \lbrace 1, \dots, n \rbrace$.
    It follows from \Cref{lemma:CrossingDirectionsMinimalGalleries} that for a minimal gallery $\gamma \subseteq \cone_v$ it holds: If $\gamma$ crosses a corresponding hyperplane $\wall_{\alpha_{v_i},k}$ with $i \in \lbrace 1, \dots, n \rbrace$  and $k \in \Z$ in $\Sigma$, then $\gamma$ crosses this hyperplane from the $\phi_{w_0}$-positive to the $\phi_{w_0}$-negative side.
    Hence, we can fold the gallery $\gamma$ with positive folding pattern $\alpha_{v_1}, \alpha_{v_2}, \cdots \alpha_{v_n}$.
\end{proof}

\begin{figure}
    \begin{subfigure}[c]{0.45\textwidth}
        \centering
    \begin{overpic}[scale=0.6]{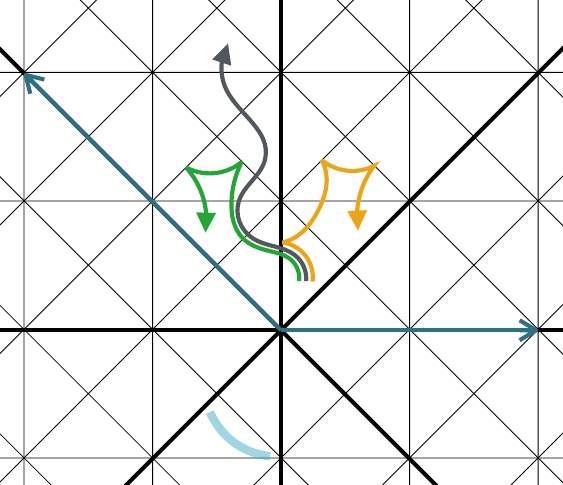}
        \put(83,29) {\tiny \textcolor{darkturquoise}{$\alpha_1$}}
        \put(11,67) {\tiny \textcolor{darkturquoise}{$\alpha_2$}}

        \put(52,33) {\tiny $e$}
        \put(43,35) {\tiny $s_1$}
        \put(35,30) {\tiny $s_1s_2$}
        \put(31,23) {\tiny $s_1s_2s_1$}
        \put(57,23) {\tiny $s_2s_1$}
        \put(58,30) {\tiny $s_2$}
        \put(43,16) {\tiny $w_0$}
        \put(50,13) {\tiny $s_2s_1s_2$}

        \put(41,70) {\footnotesize $\gamma$}
        \put(28,57) {\footnotesize $\gamma_1$}
        \put(66,57) {\footnotesize $\gamma_2$}
         %orientation
        \put(17,24) {\tiny $+$}
        \put(17,28) {\tiny $-$}
        \put(32,7) {\tiny $+$}
        \put(28,9) {\tiny $-$}
        \put(64,7) {\tiny $+$}
        \put(68,9) {\tiny $-$}
        \put(45,2) {\tiny $+$}
        \put(51,2) {\tiny $-$}
     \end{overpic}
    \end{subfigure}
    \begin{subfigure}[c]{0.45\textwidth}
        \centering
        \begin{tikzpicture}[scale=1,->, >=Stealth, shorten >=1pt, shorten <=1pt, thick]   
            %nodes     
            \node at (5, -0.7) {$e$};
            \draw[fill=black] (5,-1) circle (2pt);
            \node at (3.2, -1.7) {$s_1$};
            \draw[fill=black] (3.5,-1.7) circle (2pt);
            \node at (6.8, -1.7) {$s_2$};
            \draw[fill=black] (6.5,-1.7) circle (2pt);
            \node at (3, -3) {$s_1s_2$};
            \draw[fill=black] (3.5,-3) circle (2pt);
            \node at (7, -3) {$s_2s_1$};
            \draw[fill=black] (6.5,-3) circle (2pt);
            \node at (2.8, -4.3) {$s_1s_2s_1$};
            \draw[fill=black] (3.5,-4.3) circle (2pt);
            \node at (7.2, -4.3) {$s_2s_1s_2$};
            \draw[fill=black] (6.5,-4.3) circle (2pt);
            \node at (5, -5.3) {$s_1s_2s_1s_2=s_2s_1s_2s_1$};
            \draw[fill=black] (5,-5) circle (2pt);
            %edges
            \draw[->] (5,-1)--(3.5,-1.7);
            \draw[->] (5,-1)--(6.5,-1.7);
            \draw[->] (5,-1)--(3.5,-4.3);
            \draw[->] (5,-1)--(6.5,-4.3);
            \draw[->] (3.5,-1.7)--(3.5,-3);
            \draw[->] (3.5,-1.7)--(6.5,-3);
            \draw[->] (3.5,-1.7)--(5,-5);
            \draw[->] (3.5,-3)--(6.5,-4.3);
            \draw[->] (3.5,-3)--(3.5,-4.3);
            \draw[->] (3.5,-4.3)--(5,-5);
            \draw[->] (6.5,-1.7)--(6.5,-3);
            \draw[->] (6.5,-1.7)--(3.5,-3);
            \draw[->] (6.5,-1.7)--(6.5,-4.3);
            \draw[->] (6.5,-1.7)--(5,-5);
            \draw[->] (6.5,-4.3)--(5,-5);
            \draw[->] (6.5,-3)--(3.5,-4.3);
            %roots
            \node[color=darkturquoise] at (4.3, -1.2) {\tiny$\alpha_1$};
            \node[color=darkturquoise] at (5.7, -1.2) {\tiny$\alpha_2$};
            \node[color=darkturquoise] at (5.9, -1.85) {\tiny$\alpha_1$};
            \node[color=darkturquoise] at (4.1, -1.85) {\tiny$\alpha_2$};
            \node[color=darkturquoise] at (5.1, -1.6) {\tiny$\alpha_3$};
            \node[color=darkturquoise] at (4.85, -1.8) {\tiny$\alpha_0$};
            \node[color=darkturquoise] at (2.55, -2.4) {\tiny$2\alpha_1 + \alpha_2=\alpha_0$};
            \node[color=darkturquoise] at (7.4, -2.4) {\tiny$\alpha_3=\alpha_1+\alpha_2$};
            \node[color=darkturquoise] at (4, -2.4) {\tiny$\alpha_3$};
            \node[color=darkturquoise] at (6, -2.4) {\tiny$\alpha_0$};
            \node[color=darkturquoise] at (3.3, -3.7) {\tiny$\alpha_3$};
            \node[color=darkturquoise] at (6.7, -3.7) {\tiny$\alpha_0$};
            \node[color=darkturquoise] at (4, -4.7) {\tiny$\alpha_2$};
            \node[color=darkturquoise] at (6, -4.7) {\tiny$\alpha_1$};
            \node[color=darkturquoise] at (5.35, -3.35) {\tiny $\alpha_1$};
            \node[color=darkturquoise] at (4.7, -3.4) {\tiny $\alpha_2$};
        \end{tikzpicture}
    \end{subfigure}
    \caption{\emph{On the left:} Coxeter complex of $\atB_2$ with folded galleries. \emph{On the right:} Bruhat moment graph of $\tB_2$.}
    \label{fig:BruhatMomentGraphAndPosFoldingTypB}
\end{figure}

\begin{example}
    \Cref{fig:BruhatMomentGraphAndPosFoldingTypB} visualises \Cref{thm:positiveFoldingPatternsAndMomentGraph}.
    The grey unfolded minimal gallery $\gamma$ on the left ends in $\cone_{s_1}$.
    The galleries $\gamma_1$ and $\gamma_2$ are obtained from $\gamma$ by folding it at the fourth and fifth, or the first, fourth, and fifth panels, respectively.
    Let $\phi_{w_0}$ be the Weyl chamber orientation defined by $\partial w_0$, marked with a light blue arc.
    With respect to $\phi_{w_0}$, $\gamma_1$ is positively folded with folding pattern $(\alpha_0, \alpha_2)$.
    This corresponds to the label sequence of a directed path in $\graph_{\tB_2}$ starting in the node $s_1$.
    The gallery $\gamma_2$ is folded with folding pattern $(\alpha_1, \alpha_2, \alpha_0)$ as is obtained from $\gamma_1$ by adding the fold in the first panel.
    We read off from the moment graph $\graph_{\tB_2}$ that this folding pattern does not correspond to a label sequence of a directed path in the graph starting in $s_1$.
    The gallery $\gamma_2$ is not positively folded with respect to $\phi_{w_0}$, since the first fold is not positive.
    This also shows that adding a fold to a folded gallery before the existing folds, the positivity of these folds may be reversed.
\end{example}

With a modification of the moment graph, the results of \Cref{lemma:EdgesInMomentGraphAndW0Orientation} and \Cref{thm:positiveFoldingPatternsAndMomentGraph} can be generalized to Weyl chamber orientations $\phi_w$ with $w \neq w_0$.
Recall that in a Bruhat moment graph, directed paths correspond to minimal words in $\sW$ with the labels of the used edges giving the parallelism classes of the hyperplanes a gallery corresponding to this word is crossing.\\
Now imagine a minimal gallery in $\sW$ starting in $v \neq \id$.
Then, the end alcove of the longest minimal galleries will not be $w_0$, but the alcove $w$ maximizing $\d(v,w)$.
To represent this translation in folding patterns, we shift the labels of a Bruhat moment graph and obtain a modified moment graph:

\begin{defi}[Modified moment graph]\label{def:ModMomentGraph}
    Let $\graph_{\sW}= (V,E, \Phi^+_{\sW})$ be a Bruhat moment graph with labels $f_{{\graph}_{\sW}}(e=(u,w))= \alpha_i$.
    Now let $\varphi: V \to V$ be a bijection on the set of vertex labels and $\psi: f_{\graph_{\sW}}(E) \to f_{\graph_{\sW}^{mod}}(E)$ a bijection on the set of edge labels, such that it holds:
    \begin{align*}
        \lbrace u,w\rbrace \in E(\graph_{\sW}) &\Leftrightarrow \lbrace \varphi (u), \varphi(w)\lbrace  \in E(\graph_{\sW}^{mod}) \text{ as undirected edges},\\
        (u,w) \in E(\graph_{\sW}) &\Rightarrow (\varphi(u'), \varphi(w'))\in E(\graph_{\sW}^{mod}) \text{ with } \varphi(u')= u, \varphi(w')=w,\\
        f_{\graph_{\sW}} (e= \lbrace u,w\rbrace )= \alpha_i &\Leftrightarrow f_{\graph_{\sW}^{mod}} (\psi(e)= \lbrace \varphi(u),\varphi(w)\rbrace )= \alpha_i  \text{ as undirected edges};
    \end{align*}
    that is, the labels of the Bruhat moment graphs are rotated.
    Let $v$ be the vertex in the modified graph with no ingoing edges.
    Then we call the labeled graph $\graph_{\sW}^{mod} = (\varphi(V),E)$ with edge labels $f_{\graph_{\sW}^{mod}} (\psi(e)= (\varphi(u),\varphi(w)))= \alpha_j \in \Phi^+$ \emph{modified moment graph with minimal element $v$}.
\end{defi}

Analogous to \Cref{lemma:EdgesInMomentGraphAndW0Orientation}, we find a connection between a modified moment graph and $\phi_v$:

\begin{lemma}[Modified moment graph and $\phi_{v}$]\label{lemma:EdgesInMomentGraphAndWOrientation}
    Let $\graph_{\sW}^{mod}$ be the modified moment graph corresponding to a spherical Coxeter group $\sW$ with minimal element $w_0.v$.
    Let $(\aW,\SCox)$ be the affine Coxeter system associated with $\sW$ with Coxeter complex $\Sigma = \Sigma (\aW,\SCox)$ and let $u,w \in \sW$.
    Denote by $\partial \Sigma$ the Coxeter complex associated with $\sW$ in the boundary of $\Sigma$ and by $\partial\phi_{v}$ the Weyl chamber orientation of $\partial\Sigma$ induced by $v \in \sW$.
    Then there is a directed edge $e= (u,w) \in \graph_{\sW}^{mod}$ if and only if the Weyl chamber $u \subseteq \partial \Sigma$ is on the $\partial\phi_{v}$-negative side and $w \subseteq \partial \Sigma$ is on the $\partial\phi_{v}$-positive side of the reflection hyperplane $\wall_r$ with $w= r(u) =  u \cdot r$, $r \in R$.
\end{lemma}

\begin{proof}
    Notice that it follows from \Cref{def:ModMomentGraph} of a modified moment graph, that (undirected) incidences between vertices and their corresponding edge labels are preserved as in a Bruhat moment graph.
    Hence we have a directed edge $e= (u,w)$ or $e= (w,u)$ in a modified moment graph $\graph_{\sW}^{mod}$ with minimal element $w_0.v$ if and only if $u$ can be reflected onto $w$ with $w=r(u) = u \cdot r, \ r \in R$.
    Observe also that this edge is labeled $\alpha_i$ if and only if the reflection hyperplane $\wall_r$ is in class $\alpha_i$.
    Recall further that we have a directed edge $e= (u,w)$ in a moment graph if and only if $\d(e,u) < \d(e,w)$ (resp. $\d(u,w_0) > \d(w,w_0))$ and $w=r(u)$ for a reflection $r \in R$.
    Since in $\graph_{\sW}^{mod}$, the vertex label $w_0.v$ is mapped onto the vertex $e$ in $\graph_{\sW}$, and (undirected) incidences of vertices and edge directions are kept, we have that $\d(u,v)<\d(w,v)$ and $w = r(u)$ if and only if $e=(u,w)$ is a directed edge in $\graph_{\sW}^{mod}$ with minimal element $w_0.v$.\\
    Together with the observations above on distances in $\sW$ and the fact that $v$ is on the $\partial \phi_v$-positive side of all hyperplanes in $\partial \Sigma$, the first direction of the proof follows by the same argument as given in the proof of Lemma 4.6 (Moment graph and $\phi_{w_0}$).
    Exchange $w_0$ by $v$.\\
    To prove the converse, let $u$ be on the $\phi_{v}$-negative and $w$ be on the $\phi_{v}$-positive side of the hyperplane $\partial \wall_r$ with $w = u \cdot r$ and $r \in R$.
    Then $(u,w)$ or $(w,u)$ is an edge of $\graph_{\sW}^{mod}$.
    Since $v$ is on the $\phi_v$-positive side of all hyperplanes in $\partial \Sigma$ and $w$ is on the $\phi_{v}$-positive side of the hyperplane $\partial \wall_r$, it follows from the definition of the orientation $\partial \phi_{v}$ that $w$ and $v$ lie on the same side of the hyperplane $\partial \wall_r$ and $u$ and $v$ are separated by $\wall_r$.
    It follows that $\d(u,v)> \d(w,v)$.
    But this implies that the edge is directed towards $w$, hence $(u,w) \in E(\graph_{\sW}^{mod})$. 
\end{proof}

Having this, we are ready to generalize \Cref{thm:positiveFoldingPatternsAndMomentGraph} and, with this, state our main result:

\begin{theorem}[Positive folding patterns and modified moment graph]\label{thm:positiveFoldingPatternsAndModMomentGraph}
    Let $(\aW,\SCox)$ be an affine Coxeter system and $\sW$ the corresponding spherical Coxeter group with modified moment graph $\graph_{\sW}^{mod}= (V, E, \Phi^+_{\sW})$ with minimal element $w_0.w$.
    Denote by $\Sigma$ the affine Coxeter complex corresponding to $(\aW,\SCox)$ and let $\phi_{w}$ be the Weyl chamber orientation on $\Sigma$ with defining element $w \in \aW$.
    Let $x \in \cone_v$ and $\gamma: \cf \leadsto x\cf$ be a minimal gallery.
    Then $(\alpha_{j,i})_{i = \lbrack n \rbrack}$ is a $\phi_v$-positive folding pattern for $\gamma$ if and only if $(\alpha_{j,1}, \alpha_{j,2}, \dots, \alpha_{j,n})$ is the label index sequence of a directed path in $\graph_{\sW}^{mod}$ starting in the vertex $v$.  
\end{theorem}

\begin{proof}
    Use the same argumentation as in the proof of \Cref{thm:positiveFoldingPatternsAndMomentGraph}, substitute $\graph_{\sW}$ by $\graph_{\sW}^{mod}$, the orientation $\phi_{w_0}$ by $\phi_w$, and \Cref{lemma:EdgesInMomentGraphAndW0Orientation} by \Cref{lemma:EdgesInMomentGraphAndWOrientation}.
\end{proof}

Observe that, in general, not all $\phi_v$-positive folding patterns for a gallery $\gamma$ can be applied to it, e.g. because the gallery is too short.
Therefore, we introduce the sets $\mathcal{X}(\cone_v , \phi_w, (\alpha_{j,i})_{i \in [n]}) \subseteq \Sigma (\aW, \SCox)$ of alcoves $x\cf \in \cone_v$ that allow for the application of a $\phi_w$-positive folding pattern $(\alpha_{j,i})_{i \in [n]}$ to a minimal gallery ending in $x\cf$.
It is then a natural question to also ask for the intersection set, that is, the set of all alcoves in a common Weyl chamber, satisfying that every $\phi_w$-positive folding pattern can be realized in a folding of a gallery ending in these alcoves.

\begin{defi}
        Let $(\aW, \SCox)$ be an affine Coxeter system with its Coxeter complex $\Sigma(\aW, \SCox)$ equipped with a Weyl chamber orientation $\phi_w$.
        Let further $\cone_v$ be a Weyl chamber of $\Sigma(\aW, \SCox)$.
        We define the set $\mathcal{X} (\cone_v, \phi_w)$ of alcoves $x\cf$ of $\cone_v$, for which every $\phi_w$-positive folding pattern obtained from the modified Bruhat moment graph $\graph_{W_0}^{mod}$ starting in the vertex $v$ can be realized for at least one minimal gallery $\gamma: \cf \leadsto x\cf$ as:
        \begin{align*}
            \mathcal{X}(\cone_v, \phi_w):= \bigcap_{(\alpha_{j,i})_{i \in [n]}} \mathcal{X}(\cone_v , \phi_w, (\alpha_{j,i})_{i \in [n]}).
        \end{align*}
        Here, $\mathcal{X}(\cone_v , \phi_w,  (\alpha_{j,i})_{i \in [n]})$ denotes the set of alcoves $x\cf \in \Sigma(\aW, \SCox)$, satisfying that for the folding pattern $(\alpha_{j,i})_{i \in [n]}$ obtained from $\graph_{\sW}^{mod}$ starting in the vertex $v$ there exists a minimal gallery $\gamma: \cf \leadsto x\cf$ such that $\gamma$ can be folded following this pattern.
\end{defi}

\begin{remark}
        Observe that, restricting to folding patterns associated to directed paths in $\graph_{\sW}^{mod}$ starting in the node labeled $v$, is not a restriction when defining the sets $\mathcal{X}(\cone_v , \phi_w, (a_{j,i})_{i \in [n]})$:
        It follows from \Cref{thm:positiveFoldingPatternsAndModMomentGraph}, that folding patterns not associated with such paths cannot be applied to minimal galleries ending in $\cone_v$, hence the defined alcove sets would be empty.\\
        Observe further, that if $(a_{j,i})_{i \in [n]}$ is a $\phi_w$-positive folding pattern that continues the $\phi_w$-positive folding pattern $(a_{j,i'})_{i' \in [n']}$, then $\mathcal{X}(\cone_v , \phi_w, (a_{j,i})_{i \in [n]}) \subseteq \mathcal{X}(\cone_v , \phi_w, (a_{j,i'})_{i' \in [n']})$.
\end{remark}

Below, we present a naive construction of a set of alcoves in $\cone_v$, to which we find minimal galleries, such that every $\phi_w$-positive folding pattern obtained from $\graph_{\sW}^{mod}$ starting in node $v$ can be applied to at least one of them.
For this, we need a modified notion of a shrunken Weyl chamber as introduced in \cite[Def. 7.2.1]{Goertz_2006}:

\begin{defi}[Shrunken Weyl chamber]\label{def:ShrunkenWeylChamber}
    Denote by $\bar{\omega}_i^\vee$ the fundamental coweight corresponding to a positive simple root $\alpha_i$.
    Let $k \in \N$.
    We define the \emph{$k$-shrunken Weyl chamber $\tilde{\cone}_{w}\subset \mathcal{C}_w$} as 
    \begin{align*}
        \tilde{\cone}_{w,k} = \cone_w + k \cdot \sum_{i:\ w\alpha_i > 0} w\bar{\omega}_i^\vee .
    \end{align*}
\end{defi}

The $k$-shrunken Weyl chamber contains all the alcoves from $\cone_w$ that do not lie between the hyperplanes $H_{w\alpha, 0}$ and $H_{w\alpha,k}$ for positive simple roots $\alpha$.

\begin{proposition}[Naive subset of $\mathcal{X}(\cone_v, \phi_w)$]\label{prop:NaiveSubsetOfXcphi}
    Let $(\aW, \SCox)$ be an affine Coxeter system and its Coxeter complex $\Sigma(\aW, \SCox)$ be equipped with a Weyl chamber orientation $\phi_w$.
    Let further $\cone_v$ be a Weyl chamber of $\Sigma(\aW, \SCox)$ and $l_p$ the length of a longest directed path in $\graph_{\sW}^{mod}$ starting in the node labeled $v$.
    If $x\cf \subseteq \tilde{\cone}_{w, l_p}$, then every $\phi_w$-positive folding pattern corresponding to a directed path in $\graph_{\sW}^{mod}$ starting in $v$ can be realized in a folding of a minimal gallery $\gamma: \cf \leadsto x\cf$.
\end{proposition}

\begin{proof}
    We prove the claim for $\cone_v = \cone_f$.
    To see the proof for $\cone_v \neq \cone_f$, follow the same arguments and apply $v$ to the roots and cones.\\
    Denote by $\rho := \frac{1}{2}\sum_{\alpha \in \Phi^+}\alpha$ the half sum of positive roots.
    Observe, that a gallery $\tau: w_0\cf \leadsto \cf$ is minimal and crosses a single hyperplane of every parallelism class.
    Observe further, that $\gamma_t: \cf \leadsto t^\rho \mathbf{w_0}$ is minimal, ends in the alcove at the tip of $\tilde{\cone}_f$, and that the concatenation $\tau * \gamma_t : \cf \leadsto t^\rho \cf$ crosses at least one hyperplane of every parallelism class.
    Compare \cite[Sec. 6.2]{MilicevicSchwerThomas_2015}.
    Now let $x\cf \subseteq \tilde{\cone}_{f, l_p}$.
    Then every minimal gallery $\gamma: \cf \leadsto x\cf$ contains a minimal gallery $\gamma_{x,t}: \cf \leadsto t^{l_p \cdot \rho}\mathbf{w_0}$.
    Hence, every minimal gallery $\gamma$ crosses at least $l_p$ hyperplanes of every parallelism class.
    With this, one can construct a folded gallery for every $\phi_w$-postive folding pattern $(\alpha_{j,i})_{i \in [n]}$ obtained from $\graph_{\sW}^{mod}$ starting in the node $e$ out of $\gamma$ as follows:
    Let $\wall_1$ be the first hyperplane of class $\alpha_{j,1}$ that $\gamma$ crosses, and denote the corresponding reflection by $r_1$.
    Fold $\gamma$ at its panel contained in $\wall_1$.
    Then, the remaining part of $\gamma$ is reflected by $r_1$ and the parallelism classes of the hyperplanes $\gamma$ crosses in this reflected part may change.
    But, since $\gamma$ crosses $l_p$ hyperplanes of every parallelism class and the mapping of hyperplane classes under reflections is bijective, the remaining reflected part of $\gamma$ crosses at least $l_p -1$ hyperplanes of every parallelism class.
    Hence, find the first hyperplane $\wall_2$ of class $\alpha_{j,2}$ $\gamma$ crosses in this remaining part and fold the gallery again at its panel contained in $\wall_2$.
    With the same argument as above we see, that the remaining part of the gallery behind these two folds crosses at least $l_p-2$ hyperplanes of every parallelism class.
    Repeat the construction for the remaining elements of the folding pattern $(\alpha_{j,i})_{i \in [n]}$.
    The claim then follows by proving that $l_p \geq n$, which follows immediately from the construction of $l_p$ and \Cref{thm:positiveFoldingPatternsAndModMomentGraph}.
\end{proof}

\Cref{ex:MainThm} illustrates, that the set constructed this way is a proper subset of $\mathcal{X}(\cone_w, \phi_v)$.
We also give an example of how to construct the exact set $\mathcal{X}(\cone_w, \phi_v)$.

\begin{example}\label{ex:MainThm}
    Consider $\aW = \atA_2$ with generating set $\SCox = \lbrace s_0, s_1, s_2 \rbrace$ and positive roots $\Phi^+ = \lbrace \alpha_1, \alpha_2, \alpha_1+\alpha_2\rbrace$.
    Let the Coxeter complex $\Sigma = \Sigma (\aW, \SCox)$ be equipped with the Weyl chamber orientation $\phi_{w_0}$.
    We want to determine all $\mathbf{x} \subseteq \cone_f$, such that for every $\phi_{w_0}$-positive folding pattern we find a minimal gallery $\gamma: \cf \leadsto \mathbf{x}$ the pattern can be applied to.
    With \Cref{thm:positiveFoldingPatternsAndModMomentGraph}, we read off the longest $\phi_{w_0}$-positive folding patterns from the Bruhat moment graph: $(\alpha_1, \alpha_2, \alpha_1),\ (\alpha_1, \alpha_1+\alpha_2, \alpha_2),\ (\alpha_2, \alpha_1, \alpha_2),\ (\alpha_2, \alpha_1+\alpha_2, \alpha_1)$.
    Compare \Cref{fig:BruhatMomentGraphAndPosFolding}.
    Observe that a reflection in a fundamental hyperplane $\wall_{\alpha_i,0}$ of class $\alpha_i$ maps the other two positive roots onto each other.
    This implies, that by applying a reflection across a hyperplane of class $\alpha_i$ onto a gallery when folding it in a panel contained in this hyperplane, the hyperplane classes $\alpha_j$ of the hyperplanes crossed after this fold swap for $j\neq i$, and stay fixed for $j=i$.
    It then follows, that a gallery $\gamma: \cf \leadsto \mb{x}$ must cross hyperplanes of classes $(\alpha_1, r_{\alpha_1}(\alpha_1+\alpha_2), r_{\alpha_1}(r_{\alpha_1+\alpha_2}(\alpha_2)))=(\alpha_1, \alpha_2, \alpha_1)$ such that we can apply the folding pattern $(\alpha_1, \alpha_1+\alpha_2, \alpha_2)$.
    For the folding pattern $(\alpha_1, \alpha_2, \alpha_1)$ we derive, that the gallery must cross hyperplanes of classes $(\alpha_1, \alpha_1+\alpha_2, \alpha_2)$; for the pattern $(\alpha_2, \alpha_1, \alpha_2)$ we derive the hyperplane class sequence $(\alpha_2, \alpha_1+\alpha_2, \alpha_1)$; and for the pattern $(\alpha_2, \alpha_1+\alpha_2, \alpha_1)$ we obtain $(\alpha_2, \alpha_1, \alpha_2)$.
    Together, it follows that for alcoves $\mb{x}$ satisfying $\mb{x} \subseteq \wall_{\alpha_1,2}^+\cap \wall_{\alpha_2,2}^+\cap \wall_{\alpha_1+\alpha_2,1}^+= \wall_{\alpha_1,2}^+\cap \wall_{\alpha_2,2}^+$ all four folding patterns mentioned above can be applied to a minimal gallery $\gamma: \cf \leadsto \mb{x}$.
    \Cref{fig:ExampleNaiveSetAndExactSet} visualizes this set of alcoves.
    \begin{figure}
        \centering 
        \begin{overpic}[width=0.9\textwidth]{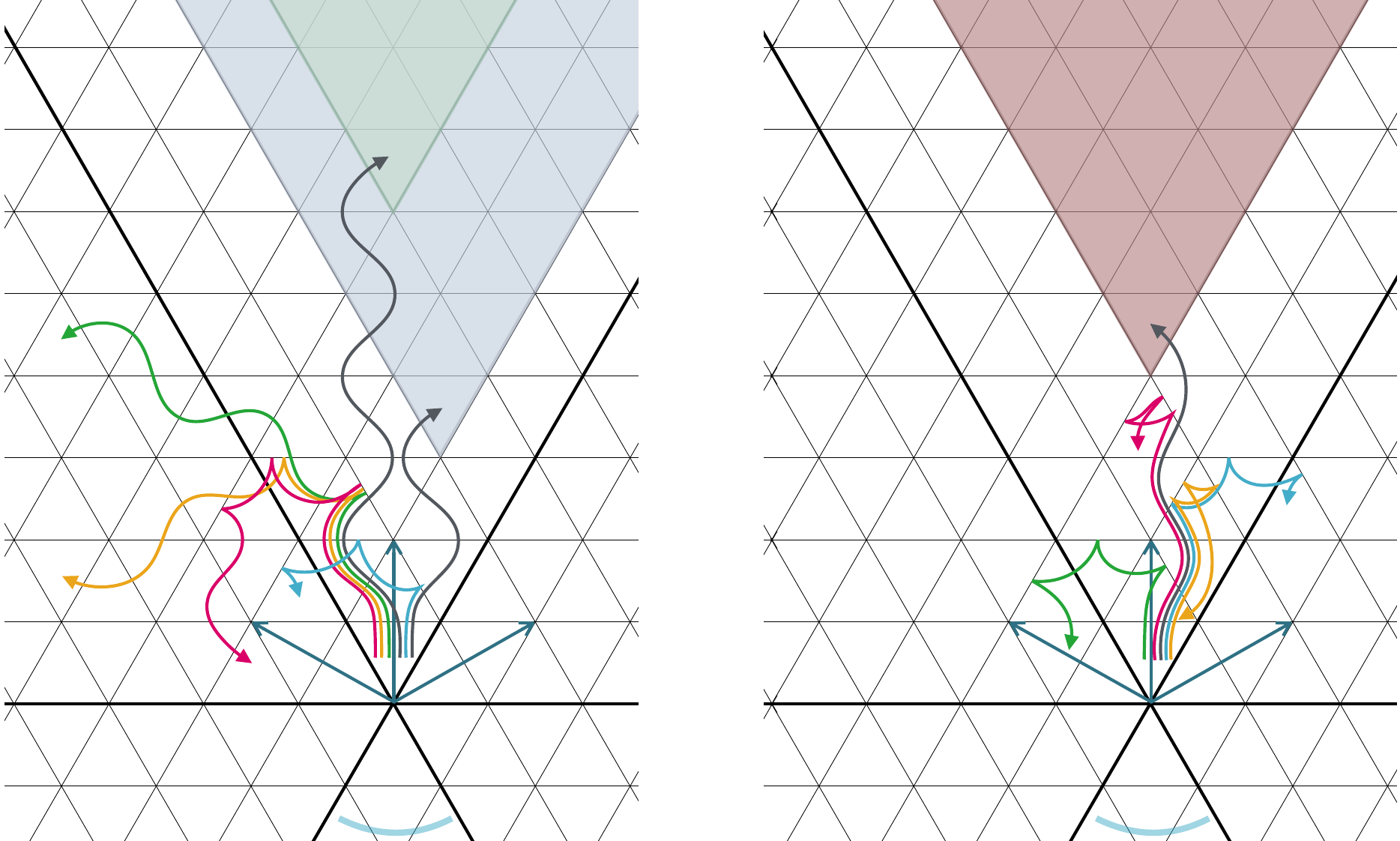}
            \put(22.5,11) {\tiny \textcolor{darkturquoise}{$\alpha_2$}}
            \put(31.5,11) {\tiny \textcolor{darkturquoise}{$\alpha_1$}}
            \put(76.5,11) {\tiny \textcolor{darkturquoise}{$\alpha_2$}}
            \put(85.5,11) {\tiny \textcolor{darkturquoise}{$\alpha_1$}}

            \put(23,48) {\footnotesize $\gamma_1$}
            \put(32,23.5) {\footnotesize $\gamma_2$}
            \put(3,34) {\footnotesize $\gamma_{1,1}$}
            \put(2,20) {\footnotesize $\gamma_{1,2}$}
            \put(12,13.5) {\footnotesize $\gamma_{1,3}$}
            \put(20,16) {\footnotesize $\gamma_{2,3}$}

            \put(27,57.5) {\footnotesize $\tilde{\cone}_{f,3}$}
            \put(83,57.5) {\footnotesize $\mathcal{X}(\cone_f, \phi_{w_0})$}
            \put(29,41) {\footnotesize $\wall_{\alpha_1,2}^+\cap \wall_{\alpha_2,1}^+$}

            %orientation
            \put(98,8) {\tiny $+$}
            \put(98,10.5) {\tiny $-$}
            \put(43.5,8) {\tiny $+$}
            \put(43.5,10.5) {\tiny $-$}
            \put(98,35) {\tiny $+$}
            \put(95,36) {\tiny $-$}
            \put(43.5,35) {\tiny $+$}
            \put(40.5,36) {\tiny $-$}
            \put(2.5,55) {\tiny $-$}
            \put(0.5,53) {\tiny $+$}
            \put(56.5,55) {\tiny $-$}
            \put(54.5,53) {\tiny $+$}
        \end{overpic}
        \caption{Visualization of \Cref{prop:NaiveSubsetOfXcphi}. The light blue arc indicates the orientation defining Weyl chamber $\partial w_0$. 
        The green cone in the picture on the left depicts the naive set constructed in \Cref{prop:NaiveSubsetOfXcphi}, and the three $i$-folded galleries $\gamma_{1,i}$ obtained from the minimal gallery $\gamma_1$, illustrate the folding procedure presented there for the folding pattern $(\alpha_1, \alpha_1+\alpha_2, \alpha_2)$.
        The blue cone depicts $\mathcal{X}(\cone_f, \phi_{w_0}, (\alpha_1, \alpha_1+\alpha_2, \alpha_2))$ as computed in \Cref{ex:MainThm}. The galleries $\gamma_2$ and $\gamma_{2,3}$ give an example of how to apply the pattern there.
        The burgundy cone in the picture on the right shows $\mathcal{X}(\cone_f, \phi_{w_0})$.
        The pink, green, yellow, and blue galleries are obtained from the grey gallery by $\phi_{w_0}$-positive folding, using all feasible folding patterns $(\alpha_1, \alpha_2, \alpha_1),\ (\alpha_1, \alpha_1+\alpha_2, \alpha_2),\ (\alpha_2, \alpha_1, \alpha_2)$, and $(\alpha_2, \alpha_1+\alpha_2, \alpha_1)$, respectively.}
        \label{fig:ExampleNaiveSetAndExactSet}
    \end{figure}
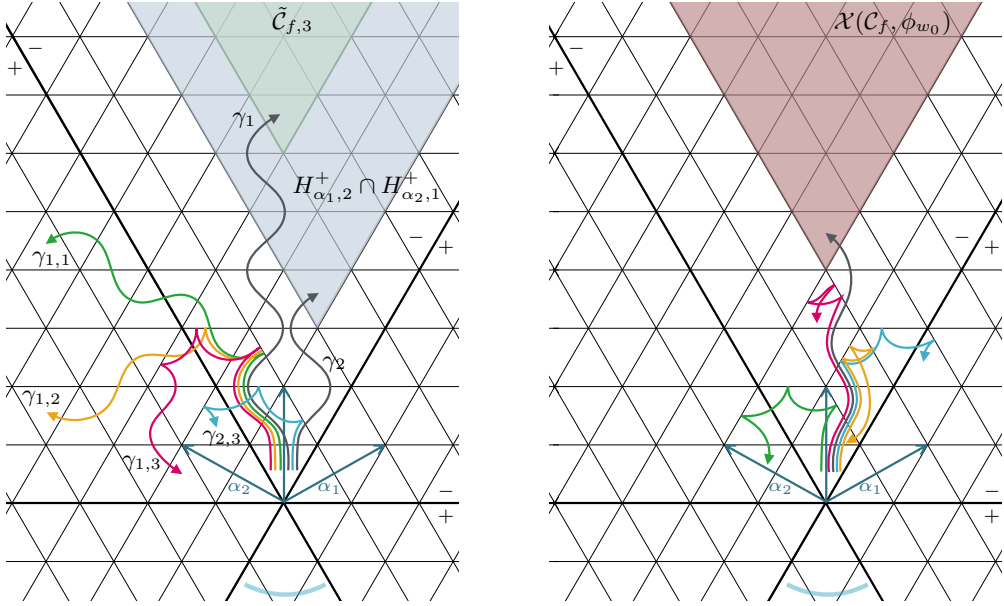
\end{example}
%
%%%%%%%%%%%%%%%%%%%%%%%%%%%%%%%%%%%%%%%%

\subsection{Moment graphs and spherical directions}

Following we present a second connection between folded galleries and Bruhat moment graphs:
Folding a gallery does not only affect the directions of crossing hyperplanes in a way that is encoded by $\graph_{\sW}$, but it also changes the spherical direction of the end alcove of a gallery.
These changes follow a structure that is encoded in an undirected version of a Bruhat moment graph, which we introduce as follows:

\begin{defi}[Undirected moment graph]\label{def:UnMomentGraph}
    Let $\graph_{\sW}= (V,E)$ with vertex set $V= \lbrace w:\ w \in \sW \rbrace$ and edges $e = (u,w) \in E$ if and only if $l(u)< l(w)$ and $w=r_{\alpha_i} (u)$ for a reflection $r_{\alpha_i} \in R$ with $\alpha_i \in \Phi^+$ a positive root in $\sW$ be the Bruhat moment graph corresponding to a spherical Coxeter system $\sW, \SCox_0$.
    Consider the undirected graph $\graph_{\sW}^{un}=(V_{un}, E_{un})$ obtained from $\graph_{\sW}= (V,E)$ by replacing every directed edge with an undirected edge between the same vertices and keeping the labels.
    Then we call this graph the \emph{undirected moment graph of $\sW$}.
\end{defi}

Observe that an undirected moment graph keeps the geometric information of a Bruhat moment graph on the adjacency of alcoves with given spherical direction in the corresponding affine Coxeter group and also gives the information on the parallelism class a panel (resp. a hyperplane) in between two such alcoves has. 
However, the ordering information on the length of the group elements in $\sW$ is lost.

\begin{theorem}[Undirected moment graph and end alcoves of folded galleries]\label{prop:UndirectedMomentGraphAndSphericalDirection}
    Let $(\aW,\SCox)$ be an affine Coxeter system with $(\sW,\SCox_0)$ the associated spherical system and let $\gamma: \cf \leadsto t^\lambda v$ be an alcove-to-alcove gallery in the affine Coxeter complex oriented by a Weyl chamber orientation $\phi_{w}$.
    Let further $(\alpha_{j,i})_{i = \lbrack n \rbrack}$ be a $\phi_w$-positive folding pattern that can be applied to $\gamma$.\\
    Consider the path in the undirected moment graph $\graph_{\sW}^{un}$ associated with $\sW$ with label sequence $(\alpha_{j,1}, \dots, \alpha_{j,n})$ starting in the vertex $v$.
    Then the end vertex of the path coincides with the spherical direction of the end alcove of every gallery $\gamma'$, obtained from $\gamma$ by folding it with the folding pattern $(\alpha_{j,i})_{i = \lbrack n \rbrack}$.
\end{theorem}

\begin{proof}
    Observe first that every vertex $v$ of an undirected moment graph $\graph_{\sW}^{un}$ satisfies $\deg (v) = \labs \Phi^+_{\sW} \rabs$ and has an adjacent edge labeled with every positive root.
    Hence every folding pattern can be realised uniquely as the label sequence of an undirected path starting in an arbitrary vertex of this graph.\\
    Observe further: Folding a gallery with end alcove $t^\lambda v$ at a panel contained in a hyperplane of class $\alpha$ is applying a reflection $r_\alpha$ with reflection hyperplane in class $\alpha$ to the end alcove.
    For this holds:
    \begin{align*}
        r_\alpha(t^\lambda v) = t^{r_{\alpha}(\lambda)}r_{\alpha,\sW} (v),
    \end{align*}
    where $r_{\alpha,\sW}$ is the reflection with reflection hyperplane in class $\alpha$ in $\sW$.
    It follows from \Cref{def:BruhatGraph} of a Bruhat graph that $r_{\alpha, \sW}(v)=: u$ is the vertex adjacent to $v$ in the undirected moment graph along the edge labeled $\alpha$.\\
    These two observations together show, that the end vertex of the path in $\graph_{\sW}^{un}$ starting in $v$ with label sequence $\alpha_{i_1}, \alpha_{i_2}, \dots, \alpha_{i_n}$ coincides with the spherical direction of the end alcove of a folded gallery $\gamma'$ obtained by folding $\gamma$ following the folding pattern.
\end{proof}

%%%%%%%%%%%%%%%%%%%%%%%%%%%%%%%%%%%%%%%%%%%%%%%%%%%%%%%%%%%%%%%%%%%%%%%%%%%%%%%%%%%%%%%%%%%%%%%%%%%%%%%%%%%%%%%%%%%%%%%%%%%%
\printbibliography
\end{document}